\documentclass[11pt, a4paper]{amsart}

\usepackage{amsmath}
\usepackage{amssymb}
\usepackage{amsthm}
\usepackage{geometry}
\usepackage{lscape}
\usepackage[utf8]{inputenc}
\usepackage{tikz}
\usepackage{tikz-3dplot}
\usetikzlibrary{calc,arrows.meta}
\usepackage{float}
\usepackage{url}
\usepackage{hyperref}
\usepackage[english]{babel}

\newtheorem{theorem}{Theorem}[section]
\newtheorem{definition}[theorem]{Definition}
\newtheorem{lemma}[theorem]{Lemma}
\newtheorem{proposition}[theorem]{Proposition}
\newtheorem{question}[theorem]{Question}

\newtheorem{claim}[theorem]{Claim}
\newtheorem{corollary}[theorem]{Corollary}
\newtheorem{remark}[theorem]{Remark}
\newtheorem{example}[theorem]{Example}

\newcommand{\dist}{\operatorname{dist}}

\newcommand{\N}{\mathbb{N}}
\newcommand{\Z}{\mathbb{Z}}
\newcommand{\A}{\mathcal{A}}
\newcommand{\B}{\mathcal{B}}

\title{The Odd Shift is not sofic in dimensions $d \ge 3$}
 \author{Nikita Gladkov}
 \address[Nikita Gladkov]{Section of Mathematics, University of Geneva, rue du Conseil-Général 7-9, 1205 Geneva.}
 \email{\texttt{nikita.gladkov@unige.ch}}
 \author{Ivan Mitrofanov}
 \address[Ivan Mitrofanov]{Section of Mathematics, University of Geneva, rue du Conseil-Général 7-9, 1205 Geneva.}
 \email{\texttt{ivan.mitrofanov@unige.ch}}
\date{\today}

\begin{document}

\begin{abstract}
    In this paper, we prove that the Odd Shift $S_{odd}$ is not sofic for dimensions $d \ge 3$. Our proof uses lower bounds from communication complexity, specifically the result of Raz and Spieker on the nondeterministic communication complexity for the connectivity problem. We also prove that $S_{odd}$ is not sofic for dimension $d=2$ under the assumption we call Mute Ocean.
\end{abstract}

\maketitle

\section{Introduction}

The classification of multidimensional shift spaces is a central theme in symbolic dynamics. While the theory of one-dimensional shifts of finite type (SFTs) and sofic shifts is well-established—characterized by forbidden blocks and regular languages, respectively—the landscape in dimensions $d \ge 2$ is significantly more complex. In high dimensions, properties such as emptiness and the extension of patterns become undecidable, and the structural limitations of sofic shifts are less understood.

A notable class of shift spaces concerns those defined by topological or geometric constraints on their configurations. Specifically, we consider the Odd Shift, denoted $S_{odd}$, defined as the set of all binary configurations on $\Z^d$ where every finite connected component of 1s has an odd cardinality. This shift serves as a natural counterpart to the ``Even Shift,'' which is known to be sofic in various contexts. In the one-dimensional case, both $S_{even}$ and $S_{odd}$ are not SFT, but are sofic. However, in higher dimensions, the geometry of connected components is complex, raising the following question.
\begin{question}\cite{H}
Is the multidimensional Odd Shift sofic? 
\end{question}

In 2011, Jullien Cassaigne announced \cite{C10} that Odd Shift is sofic in dimension $d = 2$. However, this question remained open in other dimensions~\cite{GS25, H}.  
In Theorem~\ref{thm:odd_not_sofic} we show that $S_{odd}$ is not sofic in dimensions $d \ge 3$.

\subsection{Nonsificity in dimension 3 and above}

The standard techniques for showing that a shift is non-sofic (see~\cite{P13} and the discussion below)
do not allow us to show that $S_{odd}$ is not sofic. A common approach to proving that a shift is not sofic is to show that it is not \textit{effectively closed}. 
However, $S_{odd}$ can be easily shown to be effectively closed.

In~\cite{W17}, it was noted that ``all examples known to the author of effectively closed shifts
which are not sofic were obtained by in some sense allowing elements to pack too
much important information into a small area''. We use this route, showing that connections inside the $[-n,n]^d$ cube bear more information than it is possible to transmit through its boundary. 

This way, we reduce nonsoficity to a communication complexity problem. Communication complexity is a field of theoretical computer science concerned with problems of computing a function, the inputs of which are distributed between the parties~\cite{AB09, KN97}.  We provide the overview of needed tools in section \ref{sec:preliminaries}. Unlike for computational complexity, there are often lower bounds for communication complexity of a problem. These lower bounds are sometimes useful in other models of computation. In particular, the notions from communication complexity were used to disprove soficity in papers~\cite{GJ15} and~\cite[Section~7.3.1]{C25}.

The main tool of this paper is Theorem~\ref{thm:nonsofic}. It  establishes an upper bound on the nondeterministic communication complexity of a certain problem/matrix obtained from a shift, given that the shift is sofic. 
We consider the communication problem corresponding to $S_{odd}$ in dimensions $d \ge 3$. The upper bound given by Theorem~\ref{thm:nonsofic} doesn't match the lower bound obtained by communication complexity techniques.
Specifically, we utilize the lower bound established by Raz and Spieker \cite{RS93} regarding the nondeterministic communication complexity for the connectivity problem, and, in particular, the problem of recognizing whether the union of two perfect matchings on a set of $2n$ elements is a Hamiltonian cycle.

This construction and the result of Raz and Spieker can also be used to prove the non-soficity, in dimensions $d>2$, of several other subshifts defined by constraints on finite connected components.  In particular, this applies to the One-Component Shift $S_{one}$ and One Marked Component Shift $S_{marked}$.

It is worth noting that the non-soficity of $S_{marked}$ in higher dimensions can be proved without using the rather involved theorem of Raz and Spieker; we give such a proof in Remark~\ref{rem:marked}.

We also note that Jullien Cassaigne and Pierre Guillon proved the non-soficity of One Marked Component Shift  for $d \geqslant 3$ earlier, as announced in talk \cite{CG20}.

\subsection{$S_{odd}$ is not sofic in dimension 2 under Mute Ocean assumption}

The second main result of this paper is Theorem~\ref{thm:mute_odd_not_sofic}. We show that in dimension $d = 2$ the Odd Shift is not sofic under the additional assumption we call Mute Ocean, namely that all information is passed only via 1s. In Theorem~\ref{thm:marked_not_sofic} we show the same for a related $S_{marked}$ shift.

The authors' motivation for this paper stems from problem B3 proposed at the 25th Summer Conference of the International Tournament of Towns in 2013 \cite[Problem B3]{TT13}. The problem essentially asked whether $S_{odd}$ is sofic in $d=2$. The second author was presenting the project and the first author was a participating student. The first author spent a significant amount of time trying to prove it is sofic using vector field based approaches in spirit of~\cite{CG20} with the difference that no information was transferred through 0s. So Theorem~\ref{thm:mute_odd_not_sofic} shows that these approaches were doomed to failure.

\subsection{AI Usage Disclosure}

All mathematical arguments and proofs presented in this article were developed entirely by the authors without the use of artificial intelligence. AI-based tools were used for literature review, language editing and stylistic improvement of the manuscript and did not contribute to the mathematical content, reasoning, results, or proofs.


\section{Notation and Definitions}

We adopt the standard notation for symbolic dynamics on groups, consistent with Barbieri \cite{B17}. Let $\A$ be a finite set called the \textbf{alphabet}. We consider the group $\Z^d$ acting on the set of configurations $\A^{\Z^d}$ by translation.

A \textbf{configuration} is a map $x \in \A^{\Z^d}$. For any vector $\vec{v} \in \Z^d$, the \textbf{shift map} $\sigma^{\vec{v}} : \A^{\Z^d} \to \A^{\Z^d}$ is defined by:
\[ (\sigma^{\vec{v}}(x))_{\vec{u}} = x_{\vec{u}+\vec{v}} \]

for all $\vec{u} \in \Z^d$. A \textbf{subshift} (or shift space) $X \subseteq \A^{\Z^d}$ is a subset that is closed under the product topology and invariant under the shift action.

\subsection{Ocean and Islands}

In the context of the shifts discussed in this paper, we often interpret the symbol $0 \in \A$ as the background or empty space. We refer to the set of positions holding $0$s as the \textbf{ocean}. Conversely, we refer to the connected components formed by non-zero symbols ($\A \setminus \{0\}$) as \textbf{islands}. From any cell in a connected component, one can reach any other cell in the same component by taking steps between neighboring points, that is, points that differ by exactly $1$ in precisely one coordinate while all other coordinates coincide.

\subsection{Shifts of Finite Type and Sofic Shifts}

\begin{definition}[Pattern]
	
	A \emph{pattern} is a partially defined configuration, that is, a map from a subset of $\Z^d$ (called the \emph{support} of the pattern and denoted by $supp(P)$) to the alphabet.
	
\end{definition}

\begin{definition}[Subshift of Finite Type]
	A subshift $X$ is of \textbf{finite type (SFT)} if there exists a finite set of forbidden patterns $\mathcal{F}$ such that:
	\[ X = X_\mathcal{F} := \{x \in \A^{\Z^d} : \forall \vec{v} \in \Z^d~ \forall P \in \mathcal{F}, x|_{\vec{v} + \text{supp}(P)} \neq P \} \]
\end{definition}

\begin{definition}[Sofic Shift]
	A subshift $Y \subseteq \A^{\Z^d}$ is called \textbf{sofic} if it is a factor of an SFT. That is, there exists a larger alphabet $\B$, an SFT $X \subseteq \B^{\Z^d}$, and a surjective map $\pi: \B \to \A$ (extended coordinate-wise to configurations) such that $\pi(X) = Y$.
\end{definition}

Intuitively, we can think of the alphabet $\B$ as giving the letters of $\A$ additional ``colors'' or ``hidden states.'' The forbidden patterns are defined on the colored alphabet $\B$, allowing for more complex global constraints than can be enforced on $\A$ directly.

\begin{definition}
	A \emph{nearest neighbor SFT} is an SFT that can be defined by a set of forbidden patterns $\mathcal{F}$ in which the support of each pattern consists of two neighboring cells.
\end{definition}

\begin{proposition}[{\cite[Proposition~1.2]{ABS17}, \cite[Corollary~3.43]{C25}}]\label{prop:nearest-neighbor}
	For any sofic subshift $Y$ over an alphabet $\A$ there exists a nearest neighbor SFT $X$ over an alphabet $\B$ such that $\pi(X) = Y$.
\end{proposition}

Moreover, the nearest neighbor SFT can be chosen to satisfy the stricter rule: if pairs $a_1b_1$, $a_1b_2$ and $a_2b_1$ are not prohibited, then so is $a_2b_2$. In $\mathbb{Z}^2$ this structure can also be formalized using \textbf{Marked Wang Tiles} (as done e.g. in~\cite{GJ15}). A Wang tile is a unit square with colored edges, which may be placed side-by-side only if shared edges have matching colors. A \textit{marked} Wang tile additionally carries a symbol from the alphabet $\A$ in its center. A subshift is sofic if and only if it is the set of valid configurations of central symbols formed by valid tilings of the plane.

\begin{figure}[H]
    \centering
    \begin{tikzpicture}[scale=1.5]
        \newcommand{\wangtile}[7]{
            \begin{scope}[shift={(#1,#2)}]
                \draw[thick] (0,0) rectangle (1,1);
                \node at (0.5, 0.5) {\Large \textbf{#7}};
                
                \ifnum#3=1 \fill[black] (0,1) -- (1,1) -- (0.8,0.8) -- (0.2,0.8) -- cycle; \fi
                \ifnum#4=1 \fill[black] (1,1) -- (1,0) -- (0.8,0.2) -- (0.8,0.8) -- cycle; \fi
                \ifnum#5=1 \fill[black] (0,0) -- (1,0) -- (0.8,0.2) -- (0.2,0.2) -- cycle; \fi
                \ifnum#6=1 \fill[black] (0,0) -- (0,1) -- (0.2,0.8) -- (0.2,0.2) -- cycle; \fi
                
                \draw (0,1) -- (1,1) -- (0.8,0.8) -- (0.2,0.8) -- cycle;
                \draw (1,1) -- (1,0) -- (0.8,0.2) -- (0.8,0.8) -- cycle;
                \draw (0,0) -- (1,0) -- (0.8,0.2) -- (0.2,0.2) -- cycle;
                \draw (0,0) -- (0,1) -- (0.2,0.8) -- (0.2,0.2) -- cycle;
            \end{scope}
        }

        \node at (-1.5, 2.5) {Tile for ``0'':};
        \wangtile{0}{2}{0}{0}{0}{0}{0}

        \node at (-1.5, 0.5) {Tiles for ``1'':};
        \wangtile{0}{0}{1}{0}{0}{0}{1}
        \wangtile{1.2}{0}{0}{1}{0}{0}{1}
        \wangtile{2.4}{0}{0}{0}{1}{0}{1}
        \wangtile{3.6}{0}{0}{0}{0}{1}{1}

        \wangtile{0}{-1.2}{0}{1}{1}{1}{1}
        \wangtile{1.2}{-1.2}{1}{0}{1}{1}{1}
        \wangtile{2.4}{-1.2}{1}{1}{0}{1}{1}
        \wangtile{3.6}{-1.2}{1}{1}{1}{0}{1}

    \end{tikzpicture}
    \caption{A set of Marked Wang tiles generating $S_{even}$. The tile for ``0'' has only white edges. The tiles for ``1'' have either one black edge or three black edges.}
    \label{fig:wang_tiles}
\end{figure}
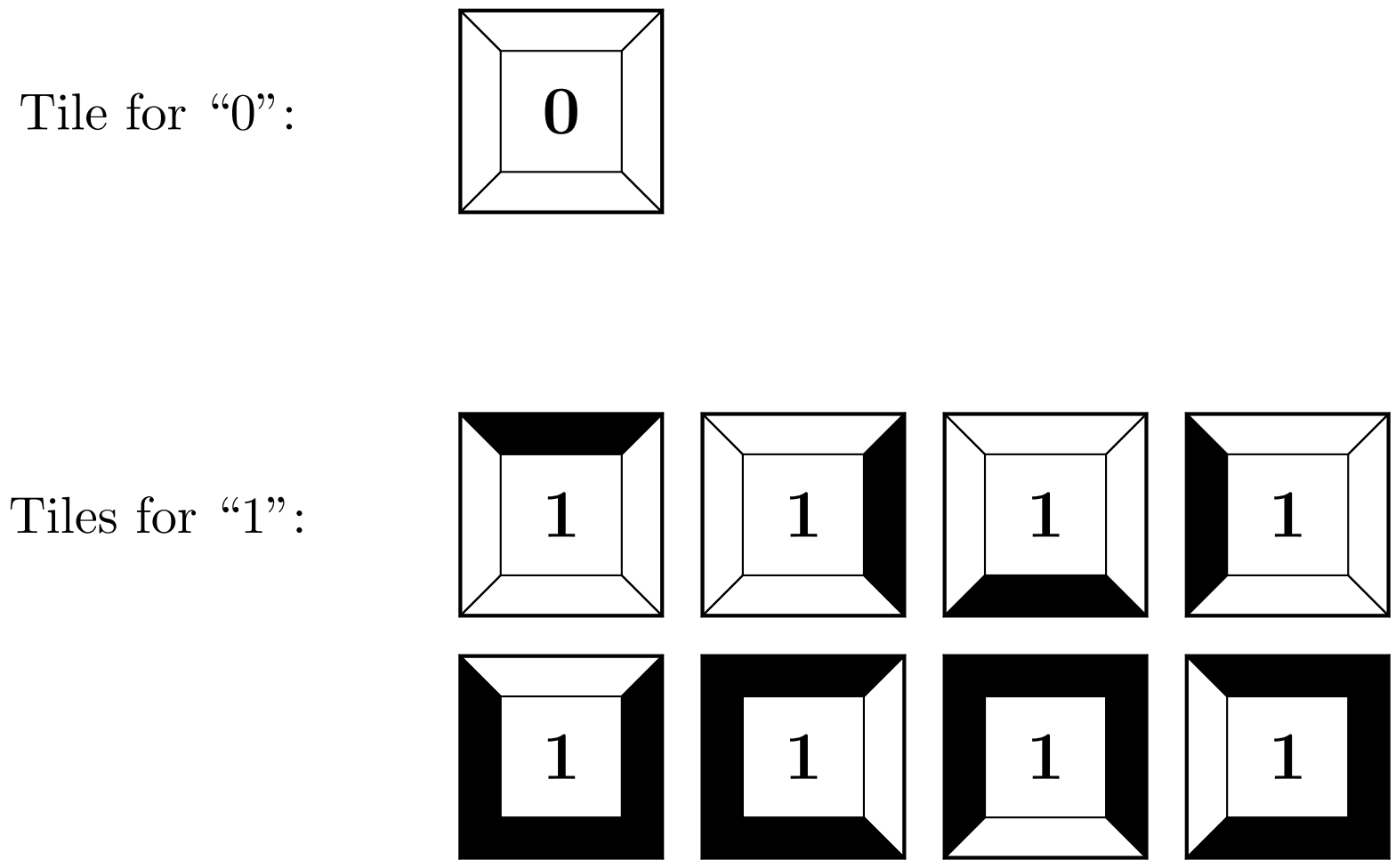

\begin{example}[$S_{even}$ in $\mathbb{Z}^2$ is Sofic]
	The ``Even Shift'' ($S_{even}$), where finite connected islands of 1s have even cardinality, is sofic. We can illustrate the valid tiles for this shift in Figure~\ref{fig:wang_tiles}. The edge colors coordinate the parity. We use a set of tiles where ``0'' (Ocean) is represented by a single tile with blank (white) edges, and ``1'' (Island) is represented by tiles with specific black/white edge configurations.
\end{example}

It is worth noting that $S_{even}$ is sofic in all dimensions and on all Cayley graphs~\cite[Prop.~1.10]{B17}. The construction relies on the well-known handshake lemma: every finite graph has an even number of odd-degree vertices. By forcing odd degrees (via edge coloring) locally, we enforce even components globally.

\subsection{The ``Mute Ocean'' Framework}

In our analysis in section \ref{sec:mute}, we utilize a simplified framework we call the \textbf{Mute Ocean}. 

\begin{definition}
    Let $A$ be an alphabet with a distinguished symbol $0_A \in A$. We say that a subshift $X \subset A^{\mathbb{Z}^d}$ admits a sofic realization with the mute ocean assumption if there exist an alphabet $B$ with a distinguished symbol $0_B \in B$, a letter projection $\pi \colon B \to A$ satisfying $\pi^{-1}(0_A)=\{0_B\}$, and a nearest-neighbor SFT $Y \subset B^{\mathbb{Z}^d}$ such that $\pi(Y)=X$.
\end{definition}

In the context of Wang tiles, this means there is exactly one marked tile with the symbol 0 in the center. As illustrated by Figure \ref{fig:wang_tiles}, $S_{even}$ is sofic even within the Mute Ocean framework.

\subsection{The One-Marked Component Shift}

Another shift of significant interest is the ``One-Marked Component'' shift, denoted here as $S_{marked}$. This shift is defined over the ternary alphabet $\A = \{0, 1, 2\}$. 
\[ S_{marked} = \{ x \in \{0, 1, 2\}^{\Z^d} \mid \forall C \in C(x), \#\{ v \in C : x_v = 2 \} = 1 \} \]
where $C(x)$ denotes the set of connected components (islands) of nonzero symbols in configuration $x$.
In other words, in every \textbf{island} (connected component of 1s and 2s), there is exactly one position labeled with the symbol ``2''.

\subsection{The Odd Shift ($S_{odd}$)}

We fix $\A = \{0, 1\}$ for this shift. We view the elements of $\Z^d$ as vertices of the standard Cayley graph induced by the generating set $\{\pm e_1, \dots, \pm e_d\}$. For a configuration $x \in \{0,1\}^{\Z^d}$, the \textbf{support} of $x$ is the set $\text{supp}(x) = \{\vec{v} \in \Z^d : x_{\vec{v}} = 1\}$. An \textbf{island} (or connected component) of $x$ is a maximal connected subgraph of the Cayley graph restricted to $\text{supp}(x)$.

We define the \textbf{Odd Shift} $S_{odd} \subset \{0,1\}^{\Z^d}$ as the set of configurations where every \textit{finite} island has odd cardinality.
\[ S_{odd} = \{ x \in \{0,1\}^{\Z^d} \mid \forall C \in C(x), |C| < \infty \implies |C| \equiv 1 \pmod 2 \} \]
where $C(x)$ denotes the set of connected components (islands) of $x$.

There is a strong structural link between these shifts. 
One can show that if $S_{marked}$ is sofic on $\Z^d$ (or any Cayley graph), then $S_{odd}$ is sofic too. Essentially, if we can sofically enforce the existence of a unique marker within each island, we can couple this with a local parity counting mechanism (similar to $S_{even}$) to enforce the odd constraint. 

\subsection{The One-Component Shift}

We define the \textbf{One-Component Shift} $S_{one} \subset \{0,1\}^{\Z^d}$ as the set of configurations with at most one finite connected component of ones.

\section{Non-Soficity of $S_{odd}$ in dimension 3 and above}

The goal of this section is to prove the following theorem.

\begin{theorem}\label{thm:odd_not_sofic}
    For $d \ge 3$, the Odd Shift $S_{odd}$ is not sofic. As a corollary, the One-Marked Component shift $S_{marked}$ and Odd-Mark Component Shift are not sofic as well.
\end{theorem}

Our proof strategy combines communication complexity lower bounds with geometric realizability. We show that if $S_{odd}$ were sofic, the number of ``border colorings'' would provide an upper bound on the nondeterministic communication complexity of a certain decision problem (see Subsection \ref{subsec:reduce}).
However, we then demonstrate that this problem requires more communication, yielding a contradiction.

\subsection{Communication Complexity Preliminaries}\label{sec:preliminaries}
We use notation consistent with the monograph~\cite{KN97}.
Let $A$ and $B$ be two sets, and let $f:A\times B \to \{0,1\}$ be a Boolean function.
A {\it rectangle} is a subset of the form $A_1\times B_1\subseteq A\times B$ for some subsets $A_1\subseteq A$ and $B_1 \subseteq B$.

\begin{definition}
	For $z\in\{0,1\}$, the {\it cover number} $C^{z}(f)$ is the minimum number $k$ such that the set
	$$f^{-1}(z) = \{(x,y)\mid f(x,y) = z\}$$
	is a union of $k$ rectangles.
\end{definition}

\begin{definition}
	Nondeterministic communication complexity of the function $f$ is $N^1(f)=\log_2 C^1(f)$.
\end{definition}

\begin{remark}
	It is often assumed that the sets $A$ and $B$ are finite. We do not make this assumption, so $N^1(f)$ may take the value $\infty$.
	However, if one of the sets is finite, then it is easy to see that $C^z(f)$ is bounded above by the cardinality of that set.
	This is exactly the situation we are interested in; see Definition~\ref{def:match}.
\end{remark}

\begin{definition}
	Associate to the function $f:A\times B\to \{0,1\}$ the matrix $M_f$ whose rows are indexed by $A$ and whose columns are indexed by $B$. The entry $(x,y)$ of $M_f$ equals $f(x,y)$. The {\it rank of the function} $rank(f)$ is the rank of the matrix $M_f$ over the field of real numbers.
\end{definition}

The reader may proceed directly to the statement of Theorem~\ref{thm:raz_spieker}, but we provide a bit more context.

\begin{definition}
	Let $f : A \times B \to \{0,1\}$. A set $S \subseteq A \times B$ is called a {\it fooling set} 	if there exists a value
	$z \in \{0,1\}$ such that
	\begin{itemize}
		\item for every $(x,y)\in S$, $f(x,y)=z$;
		\item for every two distinct pairs $(x_1,y_1)$ and $(x_2,y_2)$ in $S$, either
		$f(x_1,y_2)\neq z$ or $f(x_2,y_1)\neq z$.
	\end{itemize}
\end{definition}

It is clear that in this situation $N^z(f) \geqslant |S|$, since in any cover, two elements of $S$ cannot be covered by the same rectangle.

To help the reader internalize the definitions, we give a classical textbook example.

\begin{example}
	Let $A = B = \{0,2,\dots,2^n-1\}$, 
	$$EQ(x,y)=
	\begin{cases}
	1 \quad \text{ if } x=y\\
	0 \quad \text{ otherwise.}
	\end{cases}
	$$
	
	Then $C^1(EQ) = 2^n$, $C^0(EQ) = O(n)$, $rank(EQ) = 2^n$, $N^1(EQ) = n$.
\end{example}

\begin{proof}
	The matrix $M_{EQ}$ is the identity matrix of size $2^n\times 2^n$, hence its rank is $2^n$. Any rectangle consisting only of $1$'s has size $1\times 1$, so the set $EQ^{-1}(1)$ cannot be covered by fewer than $2^n$ rectangles (and clearly $2^n$ rectangles suffice).
	More interesting is the computation of $C^0(EQ)$.
	
	\noindent {\it Upper bound.} Let $x_i$ denote the $i$-th bit in the binary expansion of a nonnegative integer $x$.
	Then
	$$
	\{(x,y)\mid x\neq y\} = \bigcup_{i=0}^{n-1}\{(x,y)\mid x_i = 0, y_i = 1\} \cup \bigcup_{i=0}^{n-1}\{(x,y)\mid x_i = 1, y_i = 0\}.
	$$
    This gives a cover by $2n$ rectangles.
	
	\noindent {\it Lower bound.} Let $x_1,x_2\in A$. We say that a rectangle $A_1\times B_1$ {\it separates} $x_1$ and $x_2$ if the set $A_1$ contains exactly one of these two elements.
	Assume that the set $Eq^{-1}(0)=\{(x,y)\mid x\neq y\}$ is covered by fewer than $n$ rectangles. By the pigeonhole principle, there exist two distinct numbers $x_i$ and $x_j$ that are not separated by any rectangle. Some rectangle contains the pair $(x_i,x_j)$. But then the same rectangle also contains the pair $(x_j,x_j)$, a contradiction.

    \noindent {\it Exact value.} The bounds above show that $n \le C^0(EQ) \le 2n$. Note that it already determines $N_0(EQ)$ up to an \textit{additive} constant 1. However, the exact value of $C^0(EQ)$ is the smallest $m$ such that 
    $$2^n\le \binom{m}{\lfloor \frac{m}{2}\rfloor}.$$
    This result was first obtained in~\cite{CGP} (see~\cite{NC1}[Theorem 1]) as a corollary of Sperner's theorem. Asymptotically, $m$ grows like $n+O(\log n)$.
\end{proof}

\begin{remark}
    	In the literature, along with the nondeterministic communication complexity $N^1(f)$, one also studies the {\it deterministic communication complexity} $D(f)$; we refer the reader to~\cite{KN97}.
\end{remark}

\begin{remark}
There is also a game-theoretic interpretation of nondeterministic communication complexity (which we will not use):  Alice receives $x$, Bob receives $y$, and they want to compute a Boolean function $f(x,y)\in\{0,1\}$.
An oracle gives Alice and Bob a certificate $w$.
They then communicate and verify the certificate.

The protocol must satisfy:
$$
f(x,y)=1 \;\Rightarrow\; \exists w \text{ such that Alice and Bob accept},
$$
and
$$
f(x,y)=0 \;\Rightarrow\; \forall w \text{ Alice and Bob reject}.
$$
The nondeterministic communication complexity is the minimum,
over all such protocols, of the worst-case number of communicated bits
(including the certificate, under the standard convention).
\end{remark}

A large body of work is devoted to proving lower and upper bounds on the communication complexity of various functions \cite{NC1, NC2, NC3, NC4}.

Raz and Spieker proved the following theorem.

\begin{theorem}[Raz-Spieker \cite{RS93}]\label{thm:raz_spieker}
	
	Let $P_n$ be the set of all bipartite perfect matchings in the complete bipartite graph $K_{n,n}$.
	Define the function $Ham_n(x,y):P_n\times P_n\to\{0,1\}$ to be $1$ iff the union of matchings $x$ and $y$ forms a Hamiltonian cycle of length $2n$.
	
	Then
	
	\begin{enumerate}
		\item $rank(Ham_n) = \binom{2n-2}{\,n-1\,} = 2^{O(n)}$;
		\item $N^1(Ham_n) = \Omega(n\log\log n)$.
	\end{enumerate}
	
\end{theorem}

\begin{remark}
	The log-rank conjecture states that $D(f) = O((\log \text{rank}(M_f))^c)$ for some universal constant $c$, where $D(f)$ is the deterministic communication complexity. While this remains open, Raz and Spieker were the first to show that $D(Ham_N)$ and $N(Ham_N)$ are not $O(\log \text{rank}(M_{Ham_N}))$, disproving the strong form where $c=1$.
\end{remark}

The main ingredient in the proof of Theorem~\ref{thm:odd_not_sofic} is the superlinear lower bound of Raz and Spieker:
$$
\lim_{n\to \infty}\frac{N^1(Ham_n)}{n}= \infty.
$$ 

We want to use the idea that if some function is ``even more complex'' than $Ham_n$, then its communication complexity should be even larger. To make this precise, we introduce the following notion.

\begin{definition}
	Let $f_1:A_1\times B_1\to \{0,1\}$ and $f_2:A_2\times B_2\to \{0,1\}$ be two Boolean functions. We say that $f_2$ {\it reduces to} $f_1$ if there exist maps $\pi_A:A_2\to A_1$ and $\pi_B:B_2\to B_1$ such that
	$$
	f_1(\pi_A(a),\pi_B(b)) = f_2(a,b) \quad \text{ for all }(a,b)\in A_2\times B_2.
	$$
\end{definition}

\begin{lemma}\label{le:reduction}
	Let $f_2:A_2\times B_2\to \{0,1\}$ reduce to $f_1:A_1\times B_1\to \{0,1\}$. Then $N^1(f_2)\leqslant N^1(f_1)$.
\end{lemma}
\begin{proof}
	Rectangles in $A_1\times B_1$ lift naturally to rectangles in $A_2\times B_2$ via $\pi_A$ and $\pi_B$.
\end{proof}

\subsection{Soficity and communication complexity}

Let $C_n := [-n, n]^d$ be a cube in $\mathbb{Z}^d$.
Consider a configuration and two players, Alice and Bob. Alice sees only the restriction of the configuration to a cube, whereas Bob sees only the restriction to the complement of the cube. Their goal is to determine whether the configuration belongs to a given subshift. This gives rise to a natural communication problem.

\begin{definition}\label{def:match}
	Let $X\subseteq \mathcal{A}^{\Z^d}$ be a subshift, and let $n\in \N$.
	Set $A_{n}:=\mathcal{A}^{C_n}$ (all colorings of the cube), and $B_n:=\mathcal{A}^{\Z^d\setminus C_n}$ (all colorings of the complement).
	If $a\in A_n$ and $b\in B_n$, recall the notation $a\cup b \in \mathcal{A}^{\mathbb{Z}^d}$ for the configuration that equals $a$ on $C_n$ and equals $b$ outside $C_n$.
	
	Define the {\it match function} $Match_{X,n}:A_n\times B_n\to\{0,1\}$ by
	$$
	Match_{X,n}(a,b) =
	\begin{cases}
	1 \quad \text{ if } a \cup b \in X; \\
	0 \quad \text{ otherwise}.
	\end{cases}
	$$ 
\end{definition}

\begin{theorem}\label{thm:nonsofic}
Let the subshift $X\subseteq \mathcal{A}^{\Z^d}$ be sofic.
Then $N^1(Match_{X,n}) = O(n^{d-1})$.
\end{theorem}

\begin{proof}
	
	Since $X$ is sofic, there exists a subshift of finite type $Y \subseteq \mathcal{B}^{\mathbb{Z}^d}$ and a factor map $\pi: \mathcal{B} \to \mathcal{A}$ such that $\pi(Y) = X$. Moreover, by Proposition~\ref{prop:nearest-neighbor} we may assume that
	$Y$ is a nearest neighbor subshift, i.e., it is defined by a set $\mathcal{F}$ of forbidden patterns, each supported on a pair of neighboring cells.
	Let $\kappa := |\mathcal{B}|$.
	
	The boundary of $C_n$ is denoted:
	\[ \partial C_n := \{\vec{v} \in C_n : \exists i, v_i \in \{-n, n\}\} \]
	
	A \textit{border coloring} is an assignment $\omega: \partial C_n \to \mathcal{B}$. The total number of border colorings is:
	\[ |\text{BorderColorings}| = \kappa^{|\partial C_n|} \leq \kappa^{2d(2n+1)^{d-1}}. \]

	We say that a coloring of the cube $a\in A_n$ is \textit{compatible} with border coloring $\omega$ if there exists a configuration $\tilde{a} \in \mathcal{B}^{C_n}$ such that:
	\begin{enumerate}
		\item $\tilde{a}|_{\partial C_n} = \omega$ (agrees with $\omega$ on the boundary);
		\item $\tilde{a}$ contains no forbidden patterns from $\mathcal{F}$ (considering only patterns entirely within $C_n$);
		\item $\pi(\tilde{a}) = a$ (projects to $a$).
	\end{enumerate}
	
	Similarly, a coloring $b \in B_n$ is {\it compatible with $\omega$} if there exists 
	$\tilde{b} \in \mathcal{B}^{\Z^d\setminus C_n}$ 
	such that:
	\begin{enumerate}
		\item The configuration $\omega \cup \tilde{b}$ (combining the boundary with the exterior) contains no forbidden patterns from $\mathcal{F}$;
		\item $\pi(\tilde{b}) = b$.
	\end{enumerate}
	
	For each border coloring $\omega$, let $A_{\omega}$ be the set of all colorings from $A_n$ that are compatible with $\omega$, and $B_{\omega}$ be the set all colorings from $B_n$ that are compatible with $\omega$. 
	
	\begin{claim}
		\[
		Match_{X,n}^{-1}(1) = \bigcup_{\omega\in \text{BorderColorings}}A_{\omega}\times B_{\omega}
		\]
	\end{claim}
	
	\begin{proof}[Proof of Claim]
		First we show that the left-hand side is contained in the right-hand side.
		Let $a\in A_n$, $b\in B_n$, and let ${a\cup b} = x\in X.$
		Then there exists a configuration $y\in \mathcal{B}^{\mathbb{Z}^d}$ such that $\pi(y) = x$. Let $\omega$ be the restriction of $y$ to $\partial C_n$. It is immediate that $a\in A_{\omega}$ and $b\in B_{\omega}$.
		
		Now we prove the reverse inclusion. We need to check that if $a\in A_{\omega}$ and $b\in B_{\omega}$ for some $\omega\in \text{BorderColorings}$, then $a\cup b \in X$. By definition, there exist $\tilde{a} \in \mathcal{B}^{C_n}$ and $\tilde{b} \in \mathcal{B}^{\Z^d\setminus C_n}$ that both agree with $\omega$.
		We claim that the global configuration $\tilde{a} \cup \tilde{b} \in \mathcal{B}^{\mathbb{Z}^d}$ contains no forbidden patterns:
		\begin{itemize}
			\item Every forbidden pattern consists of a pair of adjacent cells, hence its support lies entirely in $C_n$ or entirely in $(\Z^d\setminus C_n) \cup \partial C_n$.
			\item Patterns entirely in $C_n$ are valid by the construction of $\tilde{a}$.
			\item Patterns entirely in $(\Z^d\setminus C_n) \cup \partial C_n$ are valid by the construction of $\tilde{b}$.
		\end{itemize}
		Therefore, $\tilde{a} \cup \tilde{b} \in Y$, which implies $\pi(\tilde{a} \cup \tilde{b}) = a \cup b \in X$. Thus $Match_{X,n}(a,b) = 1$.
	\end{proof}
	
	The theorem follows from this claim, since $\log_2(|\text{BorderColorings}|) = O(n^{d-1})$.
	
\end{proof}

The above proof of Theorem~\ref{thm:nonsofic} is fairly standard, and arguments of this kind have appeared before in work on sufficient conditions for non-soficity.
Nevertheless, we were unable to find this theorem stated explicitly in the literature.

For example, the work of Steve Kass and Kathleen Madden~\cite{KM13} gives a condition for non-soficity in terms of so-called {\it extender sets}.
Their condition admits a reformulation in the language of communication complexity.

The following notation is related to a special kind of fooling sets, the so-called {\it upper-triangular fooling sets}~\cite{KD13}.

\begin{definition}
	For a matrix $M$ let $T(A)$ be the maximum $k$ such that $A$ contains a $k\times k$ submatrix which, after permuting rows and columns, becomes lower triangular with $1$'s on the main diagonal and $0$'s above it.
	For a function $f:A\times B \to \{0,1\}$ let $T(f):=T(M_f)$.
\end{definition}

\begin{theorem}[Theorem 3.2.10 of \cite{KM13}, reformulation]\label{thm:t_nonsofic}
Let the subshift $X\subseteq \mathcal{A}^{\Z^d}$ be sofic.
Then 
	$$\log_2(T(Match_{X,n})) = O(n^{d-1}).$$
\end{theorem}

Since for any Boolean function $f$ we have $T(f) \leqslant C^1(f)$, Theorem~\ref{thm:nonsofic} implies Theorem~\ref{thm:t_nonsofic}.
Moreover, the example of the function $Ham$, the bounds of Raz and Spieker in Theorem~\ref{thm:raz_spieker}, and the obvious inequality $T(f) \leqslant rank(f)$ show that $\log_2(N^1(f))$ can grow much faster than $\log_2(T(f))$.

\subsection{Reducing the non-soficity proof to the Cord Management Problem}\label{subsec:reduce}

Assume that $S_{odd}\subset \{0,1\}^{\Z^d}$ is sofic for $d\geqslant 3$.
Consider the Boolean function $Match_{S_{odd},n}:A_n\times B_n\to \{0,1\}$.
By Theorem~\ref{thm:nonsofic},
\[
N^1(Match_{S_{odd},n}) = O(n^{d-1}).
\]

On the other hand, by Theorem~\ref{thm:raz_spieker} we have $N^1(Ham_n)=\omega(n)$.
Therefore, we obtain a contradiction with Lemma~\ref{le:reduction} once we prove the following.

\begin{claim}\label{cl:reduces}
	Let $d\geqslant 3$. There exists a constant $\alpha(d)>0$ such that for all sufficiently large $n$, the function
	$Ham_{\lfloor \alpha n^{d-1}\rfloor}$ reduces to $Match_{S_{odd},n}$.
\end{claim}

\begin{definition}[Cord in $\Z^d$]
	Let $U\subset \Z^d$. A \emph{cord} in $U$ with endpoints $x$ and $y$ is a finite sequence of lattice points
	\[
	(x=v_0, v_1,\dots, v_n=y)
	\]
	such that for every $i=1,\dots,n$ the points $v_{i-1}$ and $v_i$ are adjacent by a \emph{facet} (i.e.\ differ by $\pm 1$ in exactly one coordinate), and all vertices $v_i$ lie in $U$.
\end{definition}

\begin{definition}
	Let $U\subset \Z^d$, and let $x_1,\dots,x_{2n}$ be points on the boundary of $U$.
	Let
	\[
	P=\bigl\{\{a_1,b_1\},\dots,\{a_n,b_n\}\bigr\}
	\]
	be a perfect matching (pairing) of $\{1,\dots,2n\}$.
	We say that $P$ is \emph{realized in $U$ by $r$-disjoint cords} if there exist cords $s_1,\dots,s_n\subset U$ such that
	\begin{enumerate}
		\item for each $i$, the endpoints of $s_i$ are $x_{a_i}$ and $x_{b_i}$;
		\item the distance between any two distinct cords is at least $r$, i.e.
		\[
		\dist(s_i,s_j)\ge r\qquad \text{for all } i\neq j,
		\]
		where $\dist(A,B)=\min\{\|u-v\|:u\in A,\ v\in B\}$ ($l_\infty$-metric is used).
	\end{enumerate}
\end{definition}

\begin{theorem}\label{thm:realizability}
	Fix integers $d\ge 3$ and $r\ge 1$. Let $\Pi_n:=\{x\in C_n\mid x_1 = -n\}$ (that is, $\Pi_n$ is the discrete face of the cube $C_n$).
	There exists a constant $\alpha=\alpha(d,r)>0$ such that for every sufficiently large $n$ one can choose a subset $\mathcal{E}\subset \partial C_n$ consisting of $2T \geqslant \alpha n^{d-1}$
	in such a way that \emph{every} perfect matching on these $2T$ points can be realized inside the cube $C_n$ by $r$-disjoint cords.
\end{theorem}

We prove this theorem in the next section.
For now, we explain how Theorem~\ref{thm:realizability} implies Claim~\ref{cl:reduces}.
Set $r=10$.
Take the cube $C_n$ and choose a set $\mathcal{E}\subset \Pi_n$ consisting of $2T$ points (call them {\it exits}).
Fix a bijection between these points and the vertices of the complete bipartite graph $K_{T,T}$.
Let $P_T$ be the set of all bipartite perfect matchings in $K_{T,T}$, and let $P_1$ and $P_2$ be two copies of this set.

For each $p\in P_1$, using Theorem~\ref{thm:realizability} choose a configuration $\pi_1(p)\in A_n$ such that the cells carrying the symbol {\bf 1} form $T$ $r$-disjoint cords inside $C_n$ connecting the corresponding points of $\mathcal{E}$.
If necessary, change a few $0$-cells to $1$-cells so that each of the $1$-components (there are $T$ of them) has even cardinality.

For each point $x\in \mathcal{E}$, mark a neighboring point outside the cube and denote it by $x'$. This yields a set $\mathcal{E}'$ whose vertices also correspond to the vertices of $K_{T,T}$.
Choose one point of $\mathcal{E}$ and denote it by~$x^*$.

For each $p\in P_2$, choose a configuration in $B_n$ such that the cells carrying the symbol {\bf 1} form $r$-disjoint cords connecting the corresponding points of $\mathcal{E}'$, and such that these cords are contained in the complement $\Z^d\setminus C_n$.
Moreover, we may require that all points of the cords, except for their endpoints, lie at distance at least $2$ from the cube $C_n$.
By adding $1$-cells if needed, ensure that exactly one cord has odd cardinality (namely, the one containing the vertex~$x^*$), while all other cords have even cardinality.
Denote the resulting configuration by $\pi_2(p)\in B_n$.

Now let $p_1\in P_1$ and $p_2\in P_2$. Consider the configuration $\pi_1(p_1)\cup \pi_2(p_2)$.
Its $1$-components correspond to the connected components of the graph $p_1\cup p_2$, which is a disjoint union of cycles.
Moreover, exactly one $1$-component of $\pi_1(p_1)\cup \pi_2(p_2)$ has odd cardinality (the one containing~$x^*$).
Therefore,
$$
\pi_1(p_1)\cup \pi_2(p_2) \in S_{odd}
\ \Leftrightarrow\
p_1\cup p_2 \text{ is a Hamiltonian cycle of length } 2T.
$$

This proves Claim~\ref{cl:reduces} and, together with it, the main Theorem~\ref{thm:odd_not_sofic} (assuming that Theorem~\ref{thm:realizability} holds).
The result for $S_{marked}$ follows similarly, since the one-marked-component constraint can be reduced to the odd-parity constraint by standard constructions.

\begin{remark}
The same construction also yields two further results.
Observe also that condition
$$
p_1\cup p_2 \text{ is a Hamiltonian cycle of length } 2T
$$

is equivalent to the condition $\pi_1(p_1)\cup \pi_2(p_2) \in S_{one}$.

Therefore, the One Component Shift is also non-sofic.
\end{remark}

\begin{remark}
	This argument does not work for $d=2$, because in dimension $2$ Theorem~\ref{thm:realizability} is false: inside a square $ABCD$, any curve connecting $A$ to $C$ must intersect any curve connecting $B$ to $D$.
	Moreover, one can show that in the two-dimensional case $N^1(Match_{S_{odd},n}) = \Theta(n)$.
	
	Nevertheless, the same communication-complexity methods, combined with the result of Raz and Spieker, yield new examples of non-sofic subshifts also in dimension $2$. Consider the following set of $9$ tiles:

\medskip
	
\begin{tikzpicture}[x=1cm,y=1cm,line cap=round,line join=round]

\pgfmathsetmacro{\S}{1.6}   
\pgfmathsetmacro{\G}{0.5}   
\pgfmathsetmacro{\lw}{0.35} 
\pgfmathsetmacro{\bw}{1.2}

\pgfmathsetmacro{\xA}{0}
\pgfmathsetmacro{\xB}{\S+\G}
\pgfmathsetmacro{\xC}{2*(\S+\G)}

\pgfmathsetmacro{\yC}{0}
\pgfmathsetmacro{\yB}{\S+\G}
\pgfmathsetmacro{\yA}{2*(\S+\G)}

\tikzset{
	frame/.style={draw,line width=\lw},
	bold/.style={draw,line width=\bw}
}

\newcommand{\sq}[3]{%
	\begin{scope}[shift={(#1,#2)}]
	\draw[frame] (0,0) rectangle (\S,\S);
	#3
	\end{scope}
}

\sq{\xA}{\yA}{}

\sq{\xB}{\yA}{
	\draw[bold] (0,\S/2) -- (\S,\S/2);
}

\sq{\xC}{\yA}{
	\draw[bold] (\S/2,0) -- (\S/2,\S);
}

\sq{\xA}{\yB}{
	\draw[bold] (0,\S/2) -- (\S/2,\S/2) -- (\S/2,\S);
}

\sq{\xB}{\yB}{
	\draw[bold] (\S/2,\S) -- (\S/2,\S/2) -- (\S,\S/2);
}

\sq{\xC}{\yB}{
	\draw[bold] (\S,\S/2) -- (\S/2,\S/2) -- (\S/2,0);
}

\sq{\xA}{\yC}{
	\draw[bold] (\S/2,0) -- (\S/2,\S/2) -- (0,\S/2);
}

\sq{\xB}{\yC}{

	\draw[bold] (0,\S/2) -- (0.43*\S,\S/2);
	\draw[bold] (0.57*\S,\S/2) -- (\S,\S/2);
	\draw[bold] (\S/2,0) -- (\S/2,\S);
}

\sq{\xC}{\yC}{
	\draw[bold] (0,\S/2) -- (\S,\S/2);
	\fill (\S/2,\S/2) circle[radius=3pt];
}

\end{tikzpicture}

\medskip

The picture on the eighth tile is interpreted as indicating that the horizontal line is not broken, but rather passes slightly \emph{under} the vertical one.
The line segments passing through the centers of tiles join to form connected components.
We define the \emph{weight} of a connected component to be the number of black dots it contains (tile~9).
Consider all bi-infinite configurations in which no connected component has even weight. These configurations form a subshift, and this subshift is not sofic.
To prove non-soficity, it suffices to show that on the boundary of an $n\times n$ square one can mark $\alpha n$ points in such a way that any perfect matching between them can be drawn inside the square using this set of tiles. This is a two-dimensional analogue of Theorem~\ref{thm:realizability}.	
\end{remark}

\begin{remark}\label{rem:marked}
    We sketch an alternative proof of the non-soficity of $S_{{marked}}$ that does not use Theorem~\ref{thm:raz_spieker}.
As before, choose $c n^{d-1}$ points on the boundary of a cube of side length $n$, and consider all possible matchings of these points inside the cube (forming the set of inner configurations $A$) and outside the cube (forming the set of outer configurations $B$). On each inner matching edge, choose one marked point, while the outer matching edges contain no marked points. Then the Boolean function $Match_{S_{marked},n}$ of interest is simply the diagonal matrix of size
\[
(cn^{d-1})!! \times (cn^{d-1})!!.
\]

Hence $$
N^1(Match_{S_{marked},n}) = \log_2(cn^{d-1}!!) = \omega(n^{d-1}),
$$
and soficity of $S_{marked}$ contradicts with Theorem \ref{thm:nonsofic}.
\end{remark}

\section{Geometric Realization of Matchings}\label{sec:geometric}

In this section, we prove Theorem \ref{thm:realizability} that any perfect matching on the exits $\mathcal{E}$ can be realized by disjoint paths inside the cube $C_n$. This establishes that the domain in our communication problem is sufficiently rich to encode arbitrary matchings.

We slightly reformulate this statement, since it will be more convenient for us to use $d$ to denote not the dimension of the cube, but the dimension of the face on which the exits lie.

\begin{claim}\label{cl:realization}

For every natural number $d\geqslant 2$, there exists a constant $\alpha_d>0$
such that the following holds.
Let $C_n$ be a $(d+1)$-dimensional cube of side length $2n$.
Suppose that one of its faces is subdivided into $(2n)^d$ unit cells, and the set $\Pi_n$ of the centers of these unit cells is partitioned into pairs.

Then the points in each pair can be joined by a polygonal path
contained in $C_n$ so that the distance between any two distinct
paths is at least $\alpha_d$.
\end{claim}

We use the terminology of paper~\cite{ACG94} by Alon, Chung and Graham.

\begin{definition}

Let $G=(V,E)$ be a connected graph, where $V$ is the vertex set
and $E$ is the edge set.
Define $\operatorname{rt}(G)$, the \emph{routing number} of $G$,
to be the smallest integer $r$ such that every permutation $\pi$
in the symmetric group on $n=|V|$ symbols can be written as
$$
\pi=\tau_1\tau_2\cdots\tau_r,
$$
where each factor $\tau_i$ is a product of disjoint transpositions
of the form
$$
\tau_i=(v_1\,u_1)\cdots(v_k\,u_k),
$$
with $\{v_j,u_j\}\in E$ for every $j$.
\end{definition}

Let $K_n$ be the graph with vertex set $\Pi_n$ (the set of
$(2n)^d$ cell centers), where two vertices are adjacent if the
corresponding cells share a $(d-1)$-dimensional face.

Theorem~1 of paper~[4] states that if $G$ is a tree with $n$ vertices,
then $\operatorname{rt}(G)<3n$.
Theorem~4 of the same paper gives an upper bound for the routing
number of a Cartesian product of graphs:
$$
\operatorname{rt}(G\times G')
\leq 2\operatorname{rt}(G)+\operatorname{rt}(G').
$$

Applying these results, one can show by induction on the dimension $d$ the following

\begin{corollary}
$$
\operatorname{rt}(K_n)<6(2d-1)n.
$$
\end{corollary}

Now we are ready to prove Claim \ref{cl:realization}.

\begin{proof}
Let $\pi$ be the permutation of $\Pi_n$ that swaps the two points
in each pair. This permutation is an involution. By Corollary~4, $\pi$ can be written as
$$
\pi=\tau_r\cdots\tau_2\tau_1,
$$
where $r=6(2d-1)n$ and, for each $1\leq i\leq r$, there exists a set
$E_i$ of disjoint pairs of adjacent cells in $\Pi_n$ such that
$\tau_i$ swaps the points in each pair $(v,u)\in E_i$ and fixes
all other points.
Partition $C_n$ into $r(2n)^d$ rectangular boxes, each of dimensions
$1\times\cdots\times1\times\frac{3}{2d-1}$.

We call the last coordinate {\it vertical}. A {\it layer} is the collection of all rectangular boxes
at the same vertical level. Thus, the cube $C_n$ is partitioned
into $r$ layers, numbered from $1$ to $r$ from bottom to top.

In each layer, mark the centers of the bottom and top faces
of all the boxes. This gives two sets, each of cardinality $(2n)^d$.
The permutation $\tau_i$ defines a bijection between these sets;
join each point to its image by a line segment.

The figure below shows an example with $d=2$ and $n=2$.
The edge set $E_i$ for some $i$ is shown on the left, and the corresponding
line segments within the layer are shown on the right.

\medskip

\begin{tikzpicture}

  \begin{scope}[x=1cm,y=1cm]
 \draw[thin,step=1] (0,0) grid (4,4);
  \draw[very thick] (0,0) rectangle (4,4);

  \foreach \x/\y/\u/\v in {
    1/1/1/2,
    1/3/1/4,
    2/2/3/2,
    2/3/2/4,
    3/4/4/4,
    3/1/4/1
  } {
    \draw[blue,line width=2pt,line cap=round]
      ({\x-0.5},{\y-0.5}) -- ({\u-0.5},{\v-0.5});
  }

  \foreach \x in {1,...,4} {
    \foreach \y in {1,...,4} {
      \fill[red] ({\x-0.5},{\y-0.5}) circle[radius=1.5pt];
    }
  }
  \end{scope}

  \begin{scope}[
    xshift=6cm,
    x={(1cm,0cm)},
    y={(0.4cm,0.2cm)},
    z={(0cm,1cm)},
    gridline/.style={black,thin,opacity=0.4},
    edge/.style={black,very thick},
    hiddenedge/.style={black,thin},
    connection/.style={blue,thin,line cap=round},
    face/.style={fill=cyan,fill opacity=0.07,draw=none}
  ]
    x={(1cm,0cm)},
  y={(0.45cm,0.2cm)},
  z={(0cm,1cm)},
  gridline/.style={black,thin,opacity=0.4},
  edge/.style={black,very thick},
  hiddenedge/.style={black,thin},
  connection/.style={blue,thin,line cap=round},
  face/.style={fill=cyan,fill opacity=0.07,draw=none}
]

\foreach \h in {0,1} {
  \path[face]
    (0,0,\h) -- (4,0,\h) -- (4,4,\h) --
    (0,4,\h) -- cycle;
}

\foreach \h in {0,1} {
  \foreach \i in {1,2,3} {
    \draw[gridline] (\i,0,\h) -- (\i,4,\h);
    \draw[gridline] (0,\i,\h) -- (4,\i,\h);
  }
}

\foreach \a in {0,4} {
  \foreach \i in {1,2,3} {
    \draw[gridline] (\i,\a,0) -- (\i,\a,1);
    \draw[gridline] (\a,\i,0) -- (\a,\i,1);
  }
}

\draw[hiddenedge]
  (0,0,0) -- (0,4,0) -- (4,4,0);
\draw[hiddenedge]
  (0,4,0) -- (0,4,1);

\draw[edge]
  (0,0,1) -- (4,0,1) -- (4,4,1) --
  (0,4,1) -- cycle;
\draw[edge]
  (0,0,0) -- (4,0,0) -- (4,4,0);
\foreach \a/\b in {0/0,4/0,4/4} {
  \draw[edge] (\a,\b,0) -- (\a,\b,1);
}

\foreach \x/\y/\u/\v in {
  1/1/1/2,
  1/3/1/4,
  2/2/3/2,
  2/3/2/4,
  3/4/4/4,
  3/1/4/1
} {
  \draw[connection]
    ({\x-0.5},{\y-0.5},0) --
    ({\u-0.5},{\v-0.5},1);
  \draw[connection]
    ({\u-0.5},{\v-0.5},0) --
    ({\x-0.5},{\y-0.5},1);
}

\foreach \x/\y in {2/1,3/3,4/2,4/3} {
  \draw[connection]
    ({\x-0.5},{\y-0.5},0) --
    ({\x-0.5},{\y-0.5},1);
}

\foreach \x in {1,...,4} {
  \foreach \y in {1,...,4} {
    \foreach \h/\col in {0/red,1/orange} {
      \fill[\col]
        plot[domain=0:360,samples=40,variable=\t]
        ({\x-0.5+0.08*cos(\t)},
         {\y-0.5+0.06*sin(\t)},
         \h)
        -- cycle;
    }
  }
}

  \end{scope}

\end{tikzpicture}

\medskip

The segments in each successive layer begin where the segments
in the preceding layer end. This yields $(2n)^d$ polygonal paths,
each consisting of $r$ segments and joining a point $(x,-n)$
on the bottom face of $C_n$ to the point $(\pi(x),n)$
on the top face. 
(Here $x$ denotes all coordinates except the last.)
Since $\pi$ is an involution, the polygonal paths can be grouped
into pairs: the path from $(x,-n)$ to $(\pi(x),n)$ is paired
with the path from $(\pi(x),-n)$ to $(x,n)$.
Discard one path from each pair and add one segment to the remaining
path by joining its upper endpoint to its projection onto the bottom
face by a vertical segment.

This yields $\frac{(2n)^d}{2}$ polygonal paths, each joining
a point $(x,-n)$ to the point $(\pi(x),-n)$.

These paths intersect, but each box contains segments of at most
three distinct paths. Therefore, for some $\alpha_d>0$, the paths
can be slightly perturbed so that the distance between any two
distinct paths is at least $\alpha_d$. 
Here we use the assumption that $d>1$, since intersections
between one-dimensional curves can be easily avoided in dimension $3$.
\end{proof}

\section{Non-Soficity in the Mute Ocean Framework}\label{sec:mute}

\subsection{Case $d = 3$}

In this subsection, we establish that $S_{marked}$ and $S_{odd}$ are not sofic in dimension $d = 3$ under the Mute Ocean assumption. It obviously follows from Theorem~\ref{thm:odd_not_sofic}, however it would be useful to see the proof before going to the harder $d=2$ case. 

The argument proceeds by contradiction. Suppose that $S = S_{marked}$ is a sofic shift, and let $X \subseteq \B^{\Z^3}$ be a nearest neighbor subshift of finite type that covers $S$ via a factor map $\pi: \B \to \A$. We assume the Mute Ocean property, meaning that the background symbol $0 \in \A$ has a unique pre-image $0_\B \in \B$. Consequently, any configuration consisting of isolated one-dimensional paths of non-zero symbols (islands) surrounded by $0_\B$ is governed entirely by local constraints along the path and at the vertices.

Let $c = |\B|$ denote the size of the cover alphabet. A fundamental property of symbolic dynamics is that valid path segments of sufficient length can be extended periodically. Specifically, for any linear segment of non-zero symbols of length $L \ge c + 1$, the Pigeonhole Principle guaranties that there is a symbol that repeats twice. So, any linear segment longer than $c+1$ can be elongated by multiples of a universal pumping period 
$$N := c!$$ without introducing forbidden patterns in $X$.

\begin{figure}[H]
    \centering
    \begin{tikzpicture}[scale=0.4, 
        x={(1cm,0cm)}, 
        y={(0cm,1cm)}, 
        z={({0.5*cos(45)},{0.5*sin(45)})}] 

        \definecolor{edge1}{RGB}{220, 0, 0}      
        \definecolor{edge2}{RGB}{0, 180, 0}      
        \definecolor{edge3}{RGB}{0, 0, 220}      
        \definecolor{edge4}{RGB}{255, 140, 0}    
        \definecolor{edge5}{RGB}{148, 0, 211}    
        \definecolor{edge6}{RGB}{0, 180, 220}    

        \begin{scope}[shift={(-9, 5, 0)}]
            \node at (2, -3, 0) {\textbf{(a) Thin Cube $C_{base}$}};
            
            \coordinate (O) at (0,0,0);
            \coordinate (A) at (4,0,0);
            \coordinate (B) at (4,4,0);
            \coordinate (C) at (4,4,4);
            \coordinate (D) at (0,4,4);
            \coordinate (E) at (0,0,4);
            
            \draw[ultra thick, edge1] (O) -- (A) node[midway, below] {$e_1$};
            \draw[ultra thick, edge2] (A) -- (B) node[midway, right] {$e_2$};
            \draw[ultra thick, edge3] (B) -- (C) node[midway, right] {$e_3$};
            \draw[ultra thick, edge4] (C) -- (D) node[midway, above] {$e_4$};
            \draw[ultra thick, edge5] (D) -- (E) node[midway, left] {$e_5$};
            \draw[ultra thick, edge6] (E) -- (O) node[midway, left] {$e_6$};
            
            \fill[black] (O) circle (6pt);
            \node[below left, black] at (O) {Marked $(0,0,0)$};
        \end{scope}

        \begin{scope}[shift={(1, 0, 0)}]
            \node at (6, -3, 0) {\textbf{(b) Monster Path $P_{monster}$}};
            
            \coordinate (M0) at (0,0,0);       
            \coordinate (M1) at (8,0,0);       
            \coordinate (M2) at (8,8,0);       
            \coordinate (M3) at (8,8,8);       
            \coordinate (M4) at (4,8,8);       
            \coordinate (M5) at (4,4,8);       
            \coordinate (M6) at (4,4,4);       
            \coordinate (M7) at (12,4,4);      
            \coordinate (M8) at (12,12,4);     
            \coordinate (M9) at (12,12,12);    
            \coordinate (M10) at (0,12,12);    
            \coordinate (M11) at (0,0,12);     
            
            \draw[ultra thick, edge1] (M0) -- (M1);
            \draw[ultra thick, edge2] (M1) -- (M2);
            \draw[ultra thick, edge3] (M2) -- (M3);
            \draw[ultra thick, edge4] (M3) -- (M4);
            \draw[ultra thick, edge5] (M4) -- (M5);
            \draw[ultra thick, edge6] (M5) -- (M6);
            
            \draw[ultra thick, edge1] (M6) -- (6.4,4,4);
            \draw[ultra thick, edge1] (6.8,4,4) -- (M7);
            \draw[ultra thick, edge2] (M7) -- (M8);
            \draw[ultra thick, edge3] (M8) -- (M9);
            \draw[ultra thick, edge4] (M9) -- (M10);
            \draw[ultra thick, edge5] (M10) -- (M11);
            \draw[ultra thick, edge6] (M11) -- (M0);
            
            \fill[black] (M0) circle (6pt);
            \node[below left, black] at (M0) {Marked $(0,0,0)$};
            
            \fill[black] (M6) circle (6pt);
            \node[below left, black] at (M6) {Marked $(N,N,N)$};

            \fill[black] (4,0,0) circle (3pt);
            \fill[black] (8,4,0) circle (3pt);
            \fill[black] (8,4,4) circle (3pt);
            \fill[black] (12,8,4) circle (3pt);
            \fill[black] (8,12,12) circle (3pt);
            \fill[black] (4,12,12) circle (3pt);
            \fill[black] (0,8,12) circle (3pt);
            \fill[black] (0,4,12) circle (3pt);
            \fill[black] (0,0,4) circle (3pt);
            \fill[black] (0,0,8) circle (3pt);
            \fill[black] (8,8,4) circle (3pt);
            \fill[black] (12,12,8) circle (3pt);
            
        \end{scope}
    \end{tikzpicture}
    \caption{Color-coded correspondence between the edges of the Thin Cube and the Monster Path. The edges of the Thin Cube (a) are assigned 6 distinct colors based on their direction. The Monster Path (b) is composed of pumped versions of these edges, following the same color sequence twice. This construction forces the local environment at $(N,N,N)$ in the Monster Path to be identical to the marked environment at $(0,0,0)$ in the Thin Cube, resulting in a single component with two markers.}
    \label{fig:monster_coloring}
\end{figure}
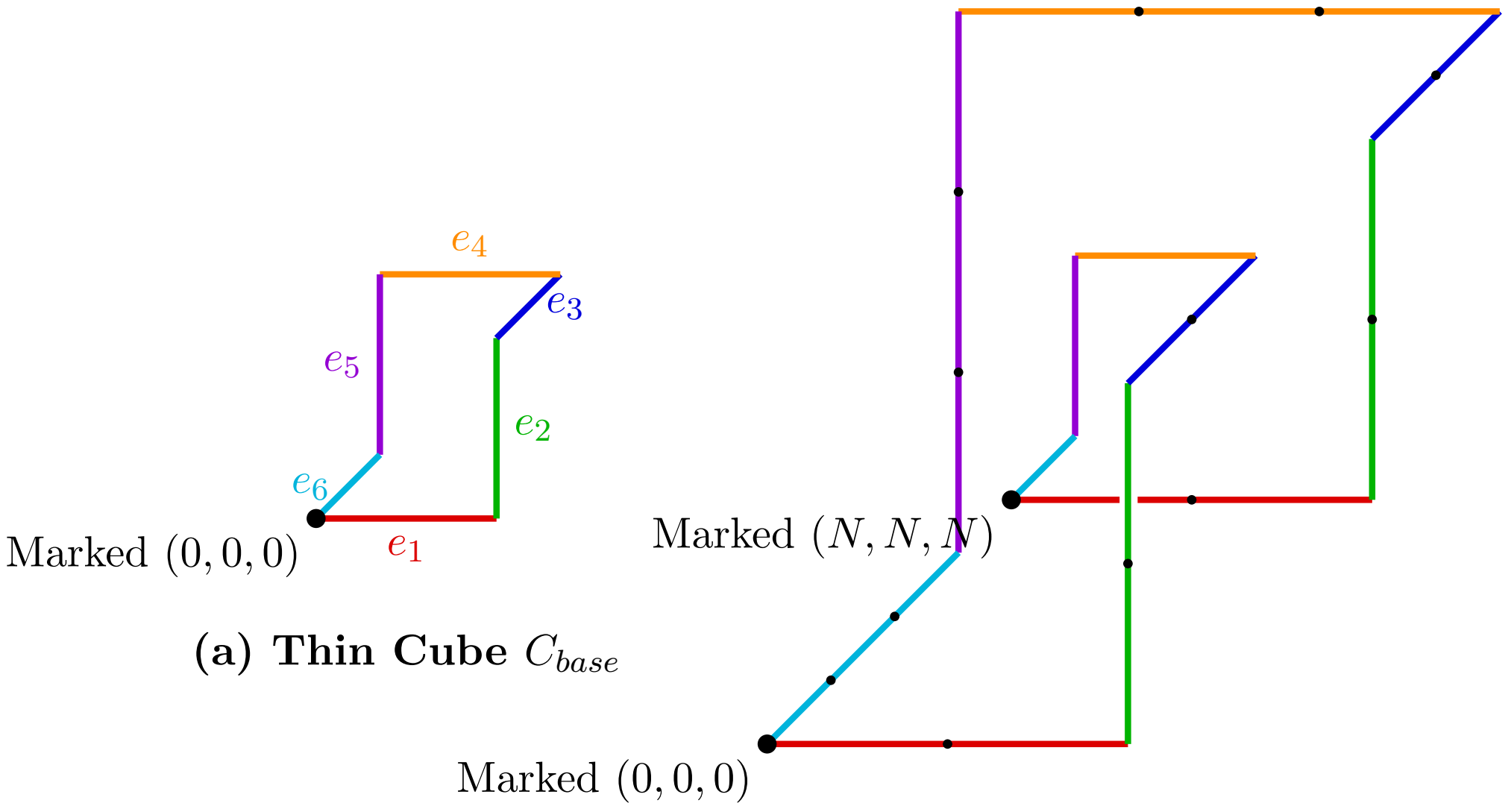

We consider a reference configuration, the ``Thin Cube'' $C_{base}$, which is a closed loop of side length $N$ visiting the vertices 
$$(0,0,0) \to (N,0,0) \to (N,N,0) \to (N,N,N) \to (0,N,N) \to (0,0,N) \to (0,0,0).$$ 
We fix the condition that the corner at $(0,0,0)$ contains the unique marker ``2''. 
That is, there is exactly one symbol $2$, exactly $6N-1$ symbols $1$, and all remaining symbols are $0$.
Since $S$ is sofic, there exists a valid pre-image $C_{colored} \in X$ such that $\pi(C_{colored}) = C_{base}$. We identify the six edges of this cube by their orientation and color them distinctively for reference, as shown in Figure \ref{fig:monster_coloring}.

We now construct a derived configuration, the ``Monster Path'' $P_{monster}\in X$, which effectively fuses two iterations of the Thin Cube geometry into a single connected component. The path is defined by the following sequence of vertices:
\[
\begin{aligned}
(0,0,0) \xrightarrow{e_1} (2N,0,0) \xrightarrow{e_2} (2N,2N,0) \xrightarrow{e_3} (2N,2N,2N) \xrightarrow{e_4} (N,2N,2N) \xrightarrow{e_5} \\
\xrightarrow{e_5} (N,N,2N) \xrightarrow{e_6} (N,N,N) \xrightarrow{e_1} (3N,N,N) \xrightarrow{e_2} (3N,3N,N) \xrightarrow{e_3} \\ 
\xrightarrow{e_3} (3N,3N,3N) \xrightarrow{e_4} (0,3N,3N) \xrightarrow{e_5} (0,0,3N) \xrightarrow{e_6} (0,0,0).
\end{aligned}
\]
Each segment of $P_{monster}$ corresponds to an edge of $C_{colored}$ elongated by a multiple of the pumping period $N$. For instance, the first segment along the $x$-axis has length $2N$, which corresponds to the first edge of $C_{colored}$ (length $N$) pumped by $N$.

Since the local validity of a configuration in $X$ is determined by permitted patterns, and the segments of $P_{monster}$ are pumped versions of the valid edges of $C_{colored}$, we conclude that $P_{monster}\in X$.
Crucially, the local geometry at the vertices is preserved. 
The neighborhoods of the vertices $(0,0,0)$ and $(N,N,N)$ in $P_{monster}$ are colored in the same way as the neighborhood of the vertex $(0,0,0)$ in $C_{colored}$. Hence, both $(0,0,0)$ and $(N,N,N)$ are marked ``2'' in $\pi(P_{monster})$.
This leads to a contradiction. The projection $\pi(P_{monster})$ constitutes a single connected component containing two distinct markers, which is forbidden in $S_{marked}$. 

An analogous argument applies to $S_{odd}$: 
let the configuration $C_{base}$ contain, instead of the cell labeled 2, a small appendage of cells labeled ``1'' in a neighborhood of the point $(0,0,0)$, so that the total number of cells in the component becomes odd.
Analogously to the case of $S_{marked}$, we construct a fusion of two such configurations, that has a preimage in $X$ but contains a single component of ``1'' of even cardinality, which is forbidden in $S_{odd}$. This contradiction shows that $S_{odd}$ can not be sofic inder the Mute Ocean condition for $d = 3$.

\subsection{Non-Soficity of $S_{marked}$ in $\mathbb{Z}^2$ under Mute Ocean assumption}

We now restrict our attention to the two-dimensional case. We show that even in $\mathbb{Z}^2$, the Mute Ocean assumption is sufficient to prevent the soficity of the One-Marked Component shift ($S_{marked}$).

\begin{theorem}\label{thm:marked_not_sofic}
    The shift $S_{marked} \subset \{0, 1, 2\}^{\mathbb{Z}^2}$ is not sofic under the Mute Ocean framework.
\end{theorem}

The idea of the proof is the same as in the three-dimensional case: construct a configuration belonging to 
$S_{marked}$ and consisting of long one-dimensional lines, find its preimage in a nearest neighbor subshift of finite type $X$, and then, using an analogue of the pumping lemma, obtain in $X$ a more complicated configuration that projects to a component with several marked vertices.

So, we begin with the assumption that $S_{marked} = \pi(X)$, where $X$ is a nearest-neighbor subshift over the alphabet $\mathcal{B}$, and $|\mathcal{B}| = c$.
As before, we denote

$$N:=c!$$

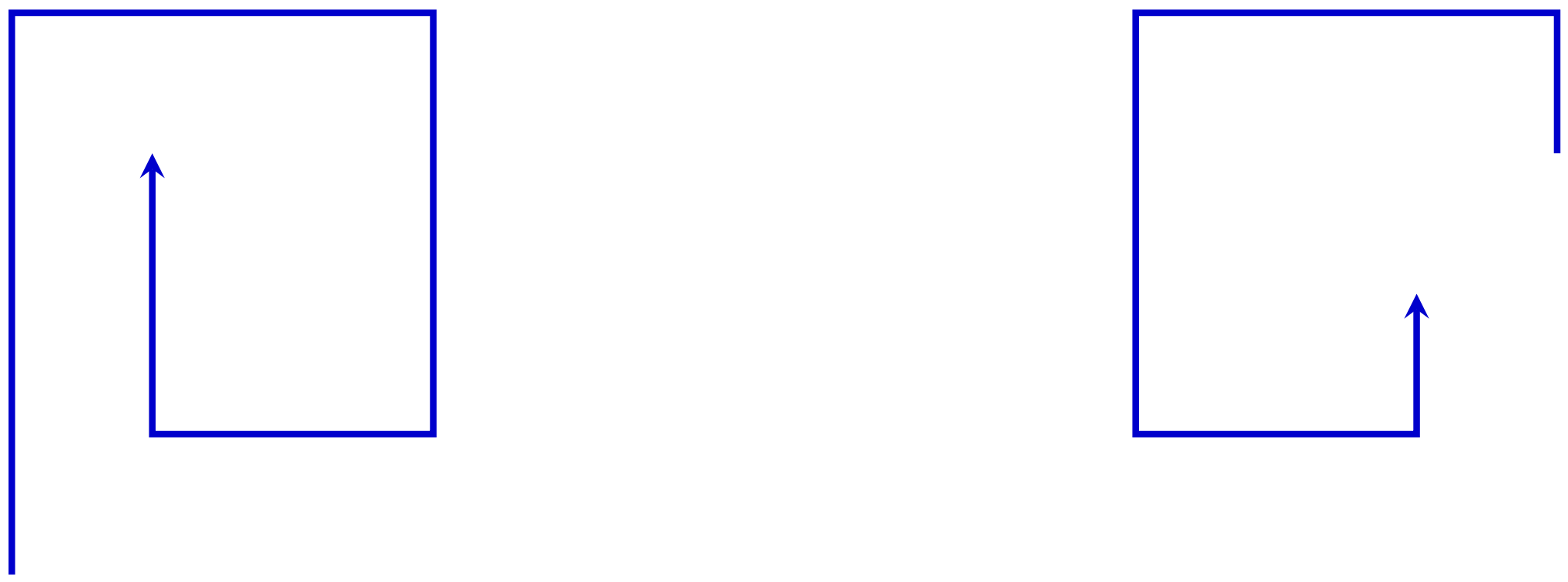
\begin{figure}[H]
    \centering
    \begin{tikzpicture}[scale=1.2]
        \tikzstyle{wire} = [ultra thick, blue!80!black, ->, >=stealth]
        
        \begin{scope}[shift={(-4,0)}]
            \node at (2, -1.5) {\textbf{(a) Clockwise Turn $T_+$}};
            \draw[wire] (0,-1) -- (0,3) -- (3,3) -- (3,0) -- (1,0) -- (1, 2);

        \end{scope}

        \begin{scope}[shift={(4,0)}]
            \node at (1, -1.5) {\textbf{(b) Counter-Clockwise Turn $T_-$}};
            \draw[wire] (3,2) -- (3,3) -- (0,3) -- (0,0) -- (2,0) -- (2, 1);
            
        \end{scope}
    \end{tikzpicture}
    \caption{Two types of turns (strictly speaking, the figures are made of cells, but we depict only the line joining the centers of the cells).}\label{fig:turns}
\end{figure}

\begin{definition}
    Fix a cell in the square grid and call it the \emph{initial cell}. A \emph{Clockwise Full Turn} is a path of a motion that starts at the inicial cell and proceeds as follows:
\begin{enumerate}
    \item move $a_1$ steps upward;
    \item move $a_2$ steps to the right;
    \item move $a_3$ steps downward;
    \item move $a_4$ steps to the left;
    \item move $a_5$ steps upward.
\end{enumerate}

We require that the path be self-avoiding and that each number $a_i$ be divisible by $N$. The final cell of the path will be called the \emph{terminal cell}.

A {\it Counter-Clockwise Full Turn} is defined analogously.

A schematic illustration of these turns is shown in Figure~\ref{fig:turns}.
\end{definition}

We will use notation $T_+$ for a Clockwise Full Turn  and $T_-$ for a Counter-Clockwise Full Turn.

These full turns can be concatenated to form closed paths; see the following figure.
We require that the initial cell of each full turn coincide with the terminal cell of the preceding full turn, and that no other self-intersections occur.
Such closed paths will henceforth be called {\it loops}.

\begin{figure}[H]
    \centering
    \begin{tikzpicture}[scale=0.8]
            \draw[ultra thick, red!80!black, ->, >=stealth] 
            (1,1)--(1,2.3) -- (2,2.3) -- (2,1) -- (1.5,1) -- (1.5,1.5);
            \draw[ultra thick, green!80!black, ->, >=stealth] 
            (1.5,1.5) -- (1.5,1.75) -- (1.25,1.75) -- (1.25, 0.5) -- (2.5, 0.5) -- (2.5, 2.5);
            \draw[ultra thick, blue!80!black, ->, >=stealth] 
            (2.5, 2.5) -- (2.5, 2.75) -- (0., 2.75) -- (0,0) -- (1,0) -- (1,1);
    \end{tikzpicture}
    \label{fig:2d_blocks}
\end{figure}

Depending on which of the three turns is chosen as the starting point, the geometry of the depicted loop can be encoded as $T_+T_-^2$, $T_-T_+T_-$, or $T_-^2T_+$.

It is easy to see that for any self-avoiding path constructed in this way, the difference between the numbers of clockwise and counterclockwise turns must be equal to $1$. The converse is also true. We give a proof similar to the one in~\cite{CR85}.

\begin{lemma}[Loop Construction]\label{lem:loop_construction}
    Let $\mathcal{S} = (s_1,s_2, \dots, s_n)$ be a sequence of symbols where each $s_i$ is either $T_+$ or $T_-$. If the total numbers of $T_+$ and $T_-$ differ exactly by one, 
    then there exists a loop in $\mathbb{Z}^2$ (a single connected component) constructed by the concatenation  of $n$ full turns, and for any $i$ the type of $i$-th turn is $s_i$.
\end{lemma}
\begin{proof}
For the moment, we ignore the condition that the segment lengths must be divisible by $N$. In fact, we work in the Euclidean plane $\mathbb{R}^2$ rather than on the square grid. Once we produce a drawing in which all segment lengths are rational, the construction of interest can then be realized on the square grid as well.

For the base case of the induction, let $n = 1$. We must construct a closed path consisting of a single turn. Such a loop is simply a rectangle, and it is clear that this can be done.

For the induction step, it suffices to show how to insert a pair $T_+T_-$ at an arbitrary position.

This is illustrated in the figure \ref{fig:insert} below.
\end{proof}

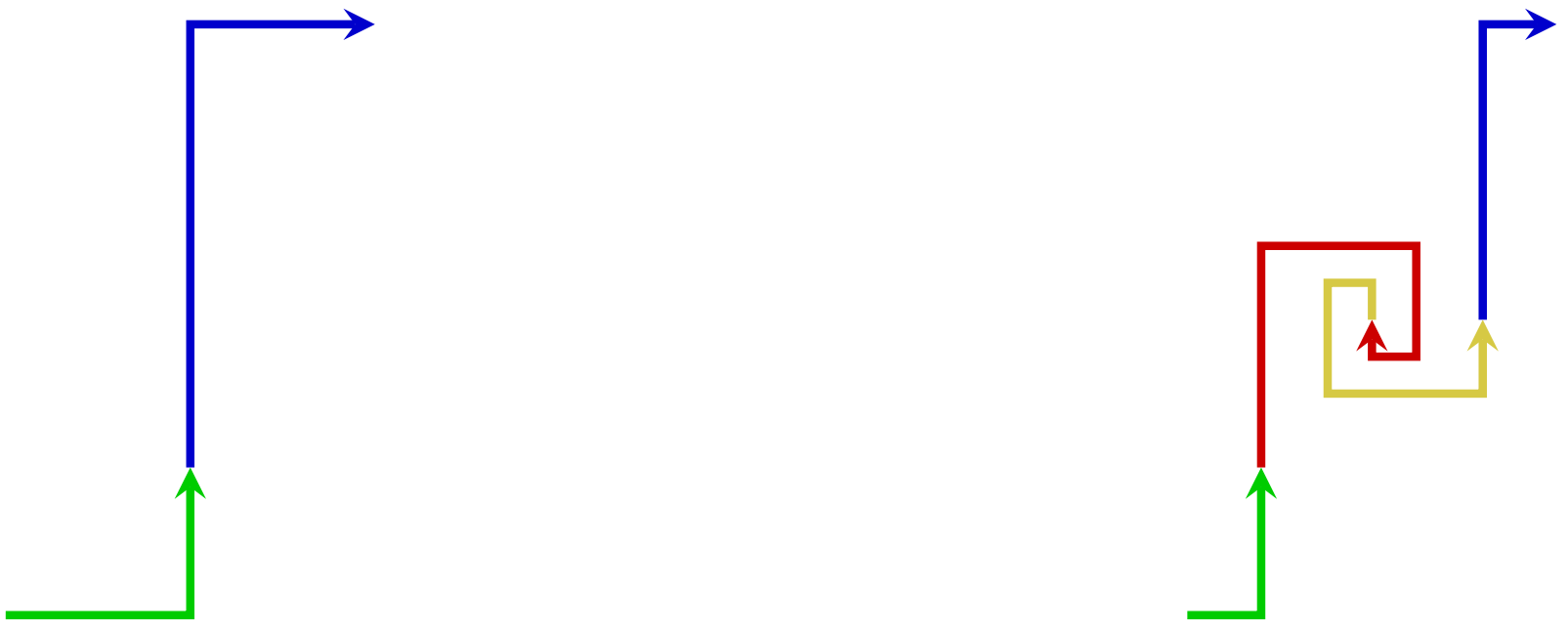
\begin{figure}[H]
    \centering
    \begin{tikzpicture}
        
        \begin{scope}[shift={(-4,0)}]
            \node at (2, -1.5) {\textbf{before inserting $T_+T_-$}};
            \draw[ultra thick, green!80!black, ->, >=stealth] 
            (-1.25,0) -- (0,0)--(0,1);
            \draw[ultra thick, blue!80!black, ->, >=stealth] 
            (0,1)--(0,4) -- (1.25,4);

        \end{scope}

        \begin{scope}[shift={(4,0)}]
            \node at (1, -1.5) {\textbf{after inserting $T_+T_-$}};
            \draw[ultra thick, green!80!black, ->, >=stealth] 
            (-1.25,0) -- (-0.75,0)--(-0.75,1);
            \draw[ultra thick, blue!80!black, ->, >=stealth] 
            (0.75,2)--(0.75,4)--(1.25,4);
            \draw[ultra thick, red!80!black, ->, >=stealth] 
            (-0.75,1)--(-0.75,2.5) -- (0.3,2.5) -- (0.3, 1.75) -- (0, 1.75) -- (0,2);
            \draw[ultra thick, yellow!80!black, ->, >=stealth] 
            (0,2) -- (0,2.25) -- (-0.3,2.25)-- (-0.3,1.5) -- (0.75,1.5) -- (0.75,2);
            
        \end{scope}
    \end{tikzpicture}
    \caption{an induction step}\label{fig:insert}
\end{figure}

\begin{definition}
Let $F \subset \mathbb{Z}^2$ be a loop obtained by concatenating full turns $T_1, \dots, T_n$.
If a configuration $x \in \mathcal{B}^{\mathbb{Z}^2}$ takes the value $0_{\mathcal{B}}$ outside $F$ and values in $\mathcal{B} \setminus \{0_{\mathcal{B}}\}$ on $F$, then we call $x$ a \emph{colored loop}.

To each such colored loop we associate its \emph{type}, which is a word of length $2n$. For each $i = 1, \dots, n$, the symbol in position $2i-1$ is either $T_+$ or $T_-$, according to the type of $T_i$, while the symbol in position $2i$ is the letter of the alphabet $\mathcal{B}$ assigned to the terminal cell of $T_i$ (equivalently, the initial cell of $T_{i+1}$).

We regard the type of a colored loop as defined up to cyclic shift.
\end{definition}

\begin{lemma}\label{le:loop_pumping}
    Let $X \subset \mathcal{B}^{\mathbb{Z}^2}$ be a nearest-neighbor subshift, and let $x \in X$ be a colored loop whose type is the cyclic word $u$.

Now let $w$ be another cyclic word such that its letters in even positions belong to $\mathcal{B}$, while its letters in odd positions are the symbols $T_+$ and $T_-$.

Assume that the following two conditions hold:
\begin{enumerate}
    \item the numbers of occurrences of $T_+$ and $T_-$ in $w$ differ by $1$;
    \item whenever a block of the form $b_1 T b_2$, where $b_1,b_2 \in \mathcal{B}$ and $T \in \{T_+,T_-\}$, appears in $w$, the same block also appears in $u$.
\end{enumerate}

Then there exists a colored loop of type $w$ that belongs to $X$.
\end{lemma}

\begin{proof}
    Let $|w| = 2n$. By Lemma~\ref{lem:loop_construction}, we construct a loop $F$ as a concatenation of $n$ full turns $T_1, \dots, T_n$ such that the type of $T_i$ is given by the $(2i-1)$-st letter of $w$.
Moreover, we may assume that every segment of these full turns is longer than any segment in the support of $y$. 

We assign the symbol $0_{\mathcal{B}}$ to every cell in $\mathbb{Z}^2 \setminus F$.
For each $i$, we place in the terminal cell of $T_i$ the letter occupying position $2i$ in $w$.

Fix $i$, and consider the full turn $T_i$. We now explain how to assign letters inside it; after that we will obtain a configuration $y$ that will be a colored loop of type $w$.

By assumption, the configuration $x$ contains a full turn $T'$ of the same type as $T_i$ and with the same pair of boundary letters. We place at the four corners of $T_i$ the same letters that appear at the corresponding corners of $T'$.

Now let $I$ be one of the five segments of $T_i$. Then there exists a segment $I'$ of $T'$ with the same orientation, the same boundary letters, and shorter length. Inside $I'$ one can find a subsegment that can be repeated, exactly as in the proof of the pumping lemma; in this way, the coloring of $I'$ can be extended to any length divisible by $N$. In particular, it can be extended to match the color $I$.

Proceeding in this manner, we fill each segment of $T_i$ with letters from the alphabet $\mathcal{B}$, and thus construct a configuration $y$. By construction, every adjacent pair of cells appearing in $y$ also appears in $x$. Therefore, $y \in X$.
\end{proof}

\begin{proof}[Proof of Theorem \ref{thm:marked_not_sofic}]
We proceed by contradiction. Assume $S_{marked}$ is sofic under the Mute Ocean assumption with a cover system $(X, \pi)$, where $X \subset \mathcal{B}^{\mathbb{Z}^2}$ is an SFT and $|\mathcal{B}| = c$. Let $L$ be a loop of type $T_+^{c+1}T_-^{c}$.
Such a loop exists because $(c+1)-c = 1$ (we apply Lemma \ref{lem:loop_construction}).
Place the symbol $2$ in the initial cell of the first full turn, assign the symbol $1$ to all other cells of $L$, and assign the symbol $0$ to all cells in $\mathbb{Z}^2 \setminus L$.
This gives a configuration $x_L \in S_{marked}$.

There exists a valid pre-image configuration $y \in X$ supported on $L$.
The configuration $y$ is a colored loop, and its type is 

$$w = (T_+b_1T_+b_2\dots T_{+}b_{c+1}T_{-}b_{c+2}\dots T_-b_{2c+1})$$

for some letters $b_1,\dots b_{2c+1}\in \mathcal{B}$. Note that $\pi(b_1) = 2$.

By the pigeonhole principle, among the letters $b_1, \dots, b_{c+1}$ there must be two that coincide. Let $b_i = b_j$, and write $j-i = k \geqslant 1$.

Consider the word obtained from $w$ by deleting the fragment $T_+b_{i+1}\dots T_+b_j$.
$$
w':= (T_+b_1T_+b_2\dots T_{+}b_iT_{+}b_{j+1}T_+\dots  b_{c+1}T_{-}b_{c+2}\dots T_-b_{2c+1}).
$$

And finally consider a word 
$$
w_{monster}:=w^kw'.
$$

In this word, the difference between the numbers of $T_+$ and $T_-$ is equal to $1$: indeed, on the block $w'$ this difference is $1-k$, while on each of the remaining $k$ blocks it is equal to $1$.

The second assumption of Lemma~\ref{le:loop_pumping} is also clearly satisfied. Therefore, there exists a colored loop $y_{monster}$ of type $w_{monster}$ belonging to the SFT $X$.

However, $\pi(b_1) = 2$, whereas $y'$ contains at least $r+1$ occurrences of the letter $b_1$. It follows that $\pi(y_{monster})$ contains more than one occurrence of the symbol $2$ and doe not belong to $S_{marked}$. 
This contradiction implies that $S_{marked}$ cannot be sofic under Mute Ocean assumption in dimension 2.
\end{proof}

\begin{figure}[H]
    \centering
    \begin{tikzpicture}[scale=0.18]
        \draw[thick, blue!50!black] 
            (52, 52) -- (52, 60) -- (59, 60) -- (59, 53) -- (54, 53) -- (54, 54) -- 
            (54, 58) -- (57, 58) -- (57, 55) -- (56, 55) -- (56, 56) -- (56, 57) -- 
            (55, 57) -- (55, 54) -- (58, 54) -- (58, 58) -- (58, 59) -- (53, 59) -- 
            (53, 52) -- (60, 52) -- (60, 60) -- (60, 61) -- (35, 61) -- (35, 34) -- 
            (62, 34) -- (62, 62) -- (62, 63) -- (17, 63) -- (17, 16) -- (64, 16) -- 
            (64, 64) -- (64, 65) -- (65, 65) -- (65, -1) -- (0, -1) -- (0, 0) -- 
            (0, 16) -- (15, 16) -- (15, 1) -- (2, 1) -- (2, 2) -- (2, 14) -- 
            (13, 14) -- (13, 3) -- (4, 3) -- (4, 4) -- (4, 12) -- (11, 12) -- 
            (11, 5) -- (6, 5) -- (6, 6) -- (6, 10) -- (9, 10) -- (9, 7) -- 
            (8, 7) -- (8, 8) -- (8, 9) -- (7, 9) -- (7, 6) -- (10, 6) -- 
            (10, 10) -- (10, 11) -- (5, 11) -- (5, 4) -- (12, 4) -- (12, 12) -- 
            (12, 13) -- (3, 13) -- (3, 2) -- (14, 2) -- (14, 14) -- (14, 15) -- 
            (1, 15) -- (1, 0) -- (16, 0) -- (16, 16) -- (16, 64) -- (63, 64) -- 
            (63, 17) -- (18, 17) -- (18, 18) -- (18, 34) -- (33, 34) -- (33, 19) -- 
            (20, 19) -- (20, 20) -- (20, 32) -- (31, 32) -- (31, 21) -- (22, 21) -- 
            (22, 22) -- (22, 30) -- (29, 30) -- (29, 23) -- (24, 23) -- (24, 24) -- 
            (24, 28) -- (27, 28) -- (27, 25) -- (26, 25) -- (26, 26) -- (26, 27) -- 
            (25, 27) -- (25, 24) -- (28, 24) -- (28, 28) -- (28, 29) -- (23, 29) -- 
            (23, 22) -- (30, 22) -- (30, 30) -- (30, 31) -- (21, 31) -- (21, 20) -- 
            (32, 20) -- (32, 32) -- (32, 33) -- (19, 33) -- (19, 18) -- (34, 18) -- 
            (34, 34) -- (34, 62) -- (61, 62) -- (61, 35) -- (36, 35) -- (36, 36) -- 
            (36, 60) -- (51, 60) -- (51, 37) -- (38, 37) -- (38, 38) -- (38, 50) -- 
            (49, 50) -- (49, 39) -- (40, 39) -- (40, 40) -- (40, 48) -- (47, 48) -- 
            (47, 41) -- (42, 41) -- (42, 42) -- (42, 46) -- (45, 46) -- (45, 43) -- 
            (44, 43) -- (44, 44) -- (44, 45) -- (43, 45) -- (43, 42) -- (46, 42) -- 
            (46, 46) -- (46, 47) -- (41, 47) -- (41, 40) -- (48, 40) -- (48, 48) -- 
            (48, 49) -- (39, 49) -- (39, 38) -- (50, 38) -- (50, 50) -- (50, 51) -- 
            (37, 51) -- (37, 36) -- (52, 36) -- (52, 52);
        
        \foreach \x/\y in {
            52/52, 52/60, 59/60, 59/53, 54/53, 54/54, 54/58, 57/58, 57/55, 56/55, 
            56/56, 56/57, 55/57, 55/54, 58/54, 58/58, 58/59, 53/59, 53/52, 60/52, 
            60/60, 60/61, 35/61, 35/34, 62/34, 62/62, 62/63, 17/63, 17/16, 64/16, 
            64/64, 64/65, 65/65, 65/-1, 0/-1, 0/0, 0/16, 15/16, 15/1, 2/1, 
            2/2, 2/14, 13/14, 13/3, 4/3, 4/4, 4/12, 11/12, 11/5, 6/5, 
            6/6, 6/10, 9/10, 9/7, 8/7, 8/8, 8/9, 7/9, 7/6, 10/6, 
            10/10, 10/11, 5/11, 5/4, 12/4, 12/12, 12/13, 3/13, 3/2, 14/2, 
            14/14, 14/15, 1/15, 1/0, 16/0, 16/16, 16/64, 63/64, 63/17, 18/17, 
            18/18, 18/34, 33/34, 33/19, 20/19, 20/20, 20/32, 31/32, 31/21, 22/21, 
            22/22, 22/30, 29/30, 29/23, 24/23, 24/24, 24/28, 27/28, 27/25, 26/25, 
            26/26, 26/27, 25/27, 25/24, 28/24, 28/28, 28/29, 23/29, 23/22, 30/22, 
            30/30, 30/31, 21/31, 21/20, 32/20, 32/32, 32/33, 19/33, 19/18, 34/18, 
            34/34, 34/62, 61/62, 61/35, 36/35, 36/36, 36/60, 51/60, 51/37, 38/37, 
            38/38, 38/50, 49/50, 49/39, 40/39, 40/40, 40/48, 47/48, 47/41, 42/41, 
            42/42, 42/46, 45/46, 45/43, 44/43, 44/44, 44/45, 43/45, 43/42, 46/42, 
            46/46, 46/47, 41/47, 41/40, 48/40, 48/48, 48/49, 39/49, 39/38, 50/38, 
            50/50, 50/51, 37/51, 37/36, 52/36
        } {
            \fill[blue!70!black] (\x, \y) circle (0.3);
        }
        
        \fill[red] (52, 52) circle (0.3);
        \node[red, above right, scale=0.8] at (52, 52){}; 
        
    \end{tikzpicture}
    \caption{The loop corresponding to the sequence $(T_+^2T_-^4T_+^5T_-^4T_+^5T_-^4T_+^5T_-^4)$. 
    The red point marks the start of the loop.}
    \label{fig:24545454}
\end{figure}
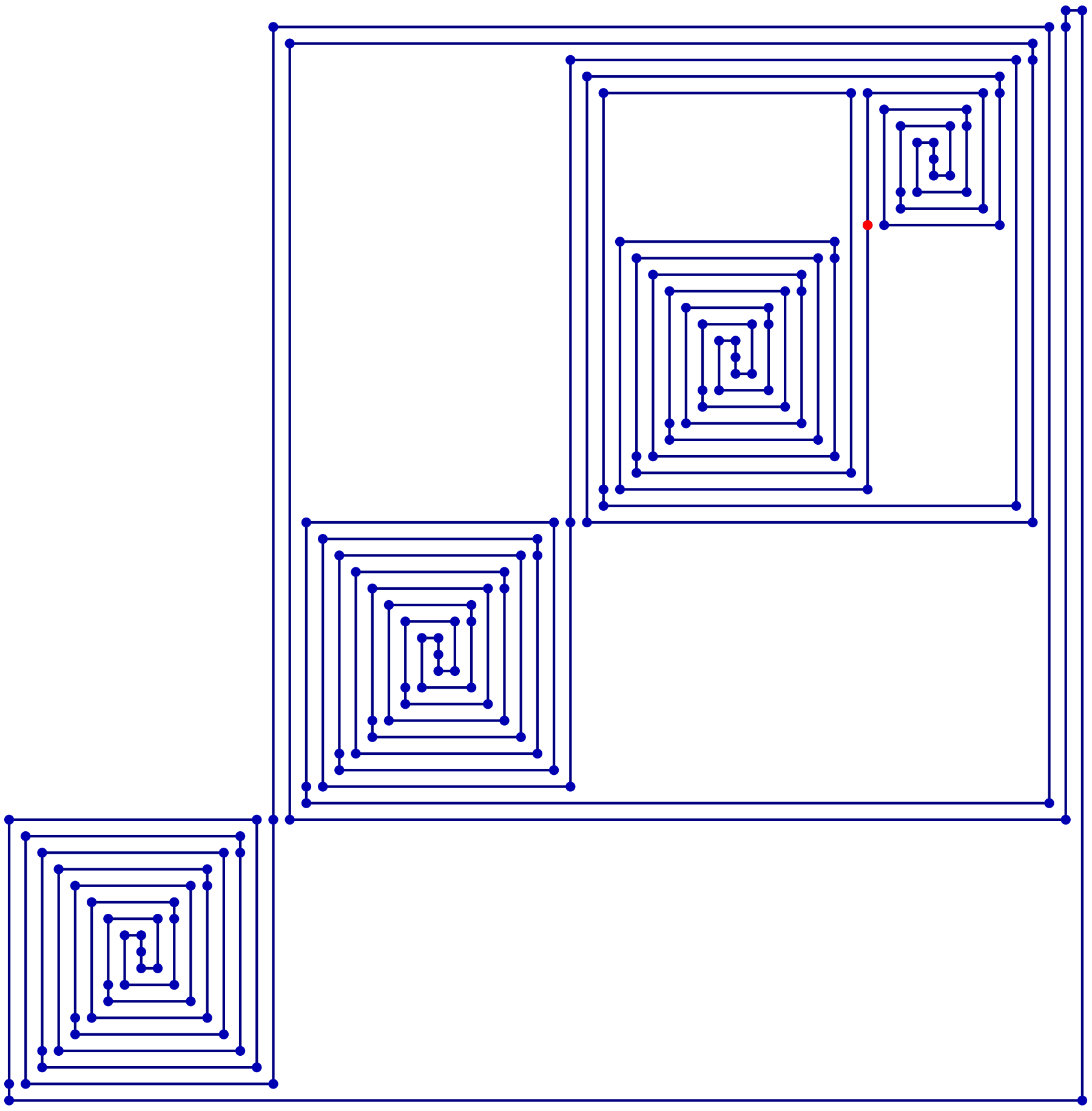

Figure~\ref{fig:24545454} illustrates the proof of the theorem and shows a possible form of $y_{\mathrm{monster}}$ in the case $c = 4$ and $k = 3$.

\subsection{Non-Soficity of $S_{odd}$ in $\mathbb{Z}^2$ under Mute Ocean assumption}

The main result of this section is

\begin{theorem}\label{thm:mute_odd_not_sofic}
    The shift $S_{odd} \subset \{0, 1\}^{\mathbb{Z}^2}$ is not sofic under the Mute Ocean framework.
\end{theorem}

The geometric constraint of parity is more elusive than the existence of a mark. In brief, we will violate it by combining two odd components into a single even component.

To prove Theorem \ref{thm:marked_not_sofic} we considered configurations whose supports were cycles, and applied some kind of pumping lemma to this cycles.
In this subsection, we will apply an analogue of the pumping lemma to more complex graphs.

\begin{definition}
Let $A$ be an alphabet that contains a letter $0_A$. 
    We call a triple $(V,E,C)$ a \emph{rectangular diagram} if the following
conditions hold:
\begin{enumerate}
    \item $V$ is a finite subset of the plane $\mathbb{R}^2$;
    \item $E$ is a set of line segments whose endpoints belong to $V$;
    \item each segment in $E$ is either horizontal or vertical;
    \item distinct segments in $E$ intersect only, if at all, at common endpoints;
    \item $C$ is a coloring of the vertices of the diagram, that is, a map
    \[
        C \colon V \to A.
    \]
\end{enumerate}
\end{definition}

We say that two rectangular diagrams are \emph{similar} if they differ only
in the lengths of their segments.

Let $K\in A^{\mathbb{Z}^2}$ be a configuration. Recall that the \emph{support} of $K$,
denoted by $\operatorname{supp} K$, is the set of cells whose symbols are
different from $0_A$.

A cell $x \in \operatorname{supp} K$ is called \emph{regular} if it has
exactly two side-neighbours in $\operatorname{supp} K$, and these neighbours
are either the upper and lower neighbours of $x$, or the right-hand and left-hand
neighbours of $x$. All other cells in $\operatorname{supp} K$ are called
\emph{singular}.

Let $U \subset \operatorname{supp} K$ be a subset containing all singular
cells. We associate with $U$ a rectangular diagram $G$ as follows. The vertices
of $G$ are the centres of the cells in $U$. For two cells $x,y \in U$, we join
the corresponding vertices by an edge if the line segment joining the centres
of $x$ and $y$ lies inside $\operatorname{supp} K$ and contains no centre of
any other cell from $U$. The coloring of the vertices is induced by the symbols
of the corresponding cells.

Assume, in addition, that the length of every segment in $G$ is divisible by
$|A|!$. We then say that the rectangular diagram $G$ is \emph{realized} by the
configuration $K$. We also say that every rectangular diagram similar to $G$ is
realized by $K$.

The \emph{geometric type} of a vertex of a rectangular diagram is the data
specifying whether there are edges incident to this vertex in each of the four
directions: up, down, right, and left. Thus, there are $16$ possible geometric
types.

The proof of the following lemma is standard, in the spirit of Lemma~\ref{le:loop_pumping}, and we omit it.

\begin{lemma}\label{le:pump}
Let $X \subset A^{\mathbb{Z}^2}$ be a nearest-neighbor SFT and let a rectangular diagram $G_1$ be realized by some configuration $K_1 \in X$.
Let $G_2$ be another rectangular diagram such that:
\begin{enumerate}
    \item for every vertex $v$ of $G_2$, there is a vertex $v'$ of $G_1$
    with the same color and the same geometric type;
    \item for every edge of $G_2$, there is an edge of $G_1$ with the same
    direction and with the same colors at its endpoints.
\end{enumerate}
Then $G_2$ is also realized by some configuration $K_2 \in X$.
Moreover, this configuration can be chosen so that the letters appearing in the
regular cells of $K_2$ form a subset of the letters appearing in the regular
cells of $K_1$.
\end{lemma}

We recall that the \textbf{Odd-Mark Component Shift} $S_{OM} \subset \{0,1,2\}^{\Z^2}$ is the set of configurations where every finite connected component of nonzero letters has odd number of letters 2.

It is straightforward to see that the soficity of $S_{odd}$ implies the soficity of $S_{OM}$
while preserving the mute ocean assumption. Therefore, for proving Theorem \ref{thm:mute_odd_not_sofic} it suffices to prove the
following result:

\begin{proposition}\label{prop:OM}
    The shift $S_{OM} \subset \{0, 1, 2\}^{\mathbb{Z}^2}$ is not sofic under the Mute Ocean framework.
\end{proposition}

As before, we start with assuming that $S_{OM} = \pi(X)$, where $X$ is a nearest-neighbor subshift over the alphabet $\mathcal{B}$, and $|\mathcal{B}| = c$. Let $N:=c!$.

\begin{figure}[hbt!]
   \includegraphics[width=0.7\linewidth]{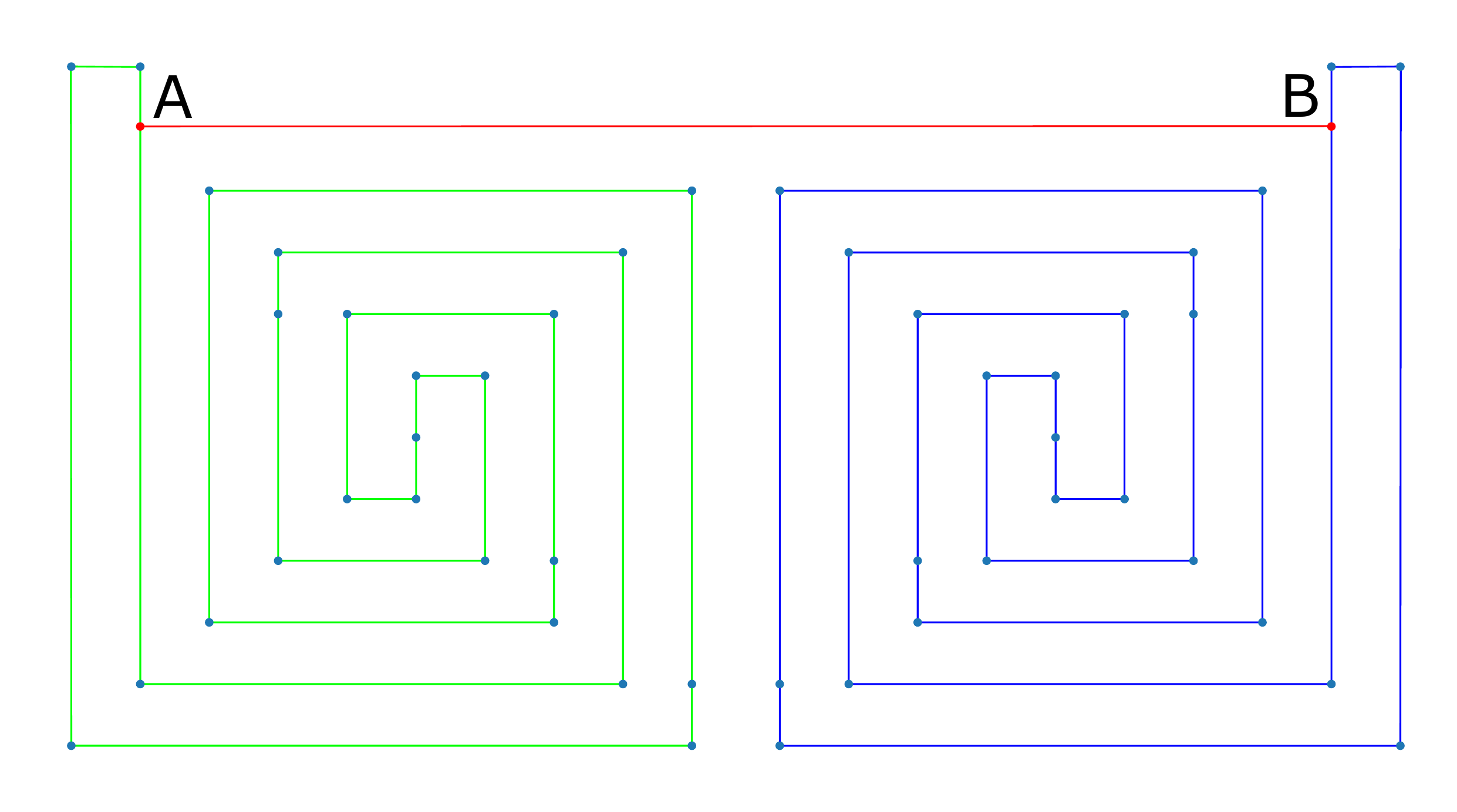}
    \caption{Rectangular diagram $F$ corresponding to $c=2$.}
    \label{fig:F}   
\end{figure}

We define a rectangular diagram $F$ with labels from the alphabet $\{0,1,2\}$
(see Figure~\ref{fig:F}). 
The figure contains two distinguished hub points, $A$ and $B$, connected by a horizontal segment $AB$. 
\begin{enumerate}
    \item At point $B$, we attach a closed loop $L$ of form $T_+^{c+1}T_-^c$. 
    \item At point $A$, we attach a symmetric loop $L'$ of form $T_+^{c}T_-^{c+1}$ (a net rotation of $-1$). 
    \item Place the letter $2$ at the vertex $A$, and the letter $1$ at all other vertices.
\end{enumerate}

The rectangular diagram $F$ is realized by some configuration $K$ containing
exactly one letter $2$ (at the position of the vertex $A$).
Clearly, the configuration $K$ is from $S_{OM}$. Hence there exists a configuration $K'$ in the SFT cover $X$ that projects to $K$.
We denote by $F'$ the rectangular diagram that corresponds to $K'$. It has the same form as $F$.

Our goal now is to construct a connected rectangular diagram $F_{\mathrm{monster}}$ such that:

\begin{enumerate}
    \item $F_{\mathrm{monster}}$ is locally arranged as $F'$ (in the sense of the conditions of Lemma~\ref{le:pump});

    \item an even number of vertices of $F_{\mathrm{monster}}$ carry a symbol that projects to $2$.
\end{enumerate}

Assume that such a diagram has already been constructed. By Lemma~\ref{le:pump}, the diagram $F_{\mathrm{monster}}$ is realized by some configuration $K_{\mathrm{monster}} \in X$ such
that all regular cells of its support project under $\pi$ to $1$. Therefore, the
support of the configuration $\pi(K_{\mathrm{monster}})$ is a connected set containing an even
number of symbols $2$, and hence $\pi(K_{\mathrm{monster}}) \notin S_{OD}$. This gives a
contradiction.

Before constructing $F_{monster}$, let us look at the diagram $F'$. 
We apply the Pigeonhole Principle to the ``necks'' of the loops. The loop $L$ at $B$ consists of $c+1$ blocks of $T^+$. Let the colors of the $T^+$ blocks at the start of the loop be $b_1, b_2, \dots, b_{c+1}$. Since there are more blocks than colors in $\mathcal{B}$, two colors in this sequence repeat.  Let $k$ be the period of this repetition. Analogously, $k'$ is the period of a repetition in $L'$.

\begin{figure}[hbt]
   \includegraphics[width=0.99\linewidth]{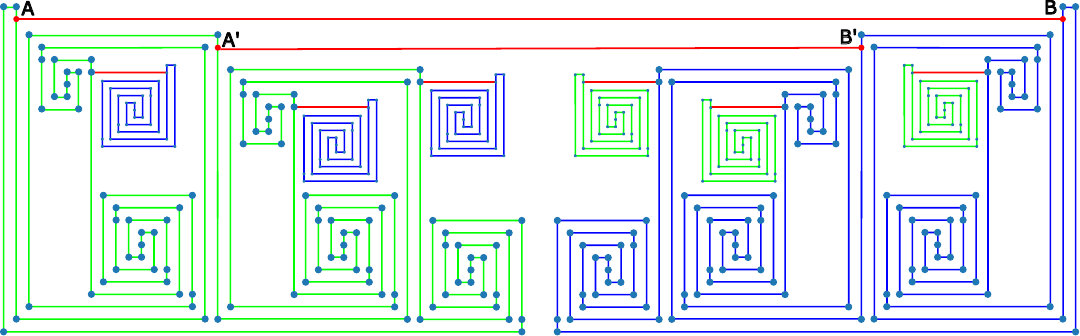}
    \caption{Monster Figure $F_{monster}$ for $c=2$, $k=k'=2$.}
    \label{fig:monster_z2}   
\end{figure}

We now construct the {Monster Figure} $F_{monster}$, which fuses two copies of $F'$ into a single component. The geometry is depicted in Figure \ref{fig:monster_z2}. The figure is defined by four hub points connected by paths: $A, B$ (top) and $A', B'$ (bottom).

\begin{enumerate}
    \item \textbf{Horizontal Bars:} $A$ is connected to $B$, and $A'$ is connected to $B'$ by horizontal segments.
    
    \item \textbf{Right Path ($B \to B'$):} This path connects the top-right hub to the bottom-right hub. It behaves like a monster-loop from the proof of Theorem \ref{thm:marked_not_sofic}. It has form:
    \[ \underbrace{T_+^{c+1} T_-^{c} \dots T_+^{c+1} T_-^{c}}_{k \text{ times}} T_+^{c+1-k} T_-^{c}, \]
    where $k$ is a valid pumping length in the right loop of $F'$.
    This effectively repeats the loop pattern $k$ times. 
    All $k$ turning points at which the sequence of blocks $T_-$ changes to $T_+$ will have the same color as the vertex $B$ of figure $F'$. 
    We also want all these $k$ points to have the same geometric type as $B$. For this reason we attach at each such a point a small copy of the \textit{left part} of the original figure $F'$ (the horizontal line to a point representing hub $A$ and attached loop $L'$).

    \item \textbf{Right Return Path ($B' \to B$):} To close the loop, we connect $B'$ back to $B$. This path has turning number 0, constructed via the sequence:
    \[ \underbrace{T_+^{c+1} T_-^{c} \dots T_+^{c+1} T_-^{c}}_{k-1 \text{ times}}T_+^{c+1-k}T_-^{c} .\]

    As before, at each switch from $T^-$ to $T^+$ we add a small copy of the left part of $F'$ ($k-1$ of them).

    \item \textbf{Left Path ($A \to A'$):} This side is constructed symmetrically to the right path, using the monster-loop logic. The unintuitive difference is that instead of $k$ we use $k'$ -- a valid pumping length in loop $L'$ in the coloring $x_F$. We can encode this monster-loop in our notation as
    
    \[ 
    \underbrace{T_-^{c+1} T_+^{c} \dots T_-^{c+1} T_+^{c}}_{k' \text{ times}} 
    T_-^{c+1-k'} T_+^{c}. 
    \]

    At every switch from $T^+$ to $T^-$ (there are $k'$ of them) we add a small copy of the {\it right part} of $F$, these $k'$ points have the same local structure as the point $A$ of $F'$.

    \item \textbf{Left Return Path ($A' \to A$):}
    The left return path is constructed symmetrically to the right return path. We can encode this loop in our notation as
    \[ \underbrace{T_-^{c+1}T_+^{c} \dots T_-^{c+1}T_+^{c}}_{k'-1 \text{ times}} T_-^{c+1-k'}T_+^{c}.\]
    
    At every switch from $T^+$ to $T^-$ we add a small copy of the right part of $F'$.
\end{enumerate}

Strictly speaking, we should explain why all these paths can be drawn without introducing self-intersections. This is intuitively clear from Figure~\ref{fig:monster_z2}. A formal proof can be easily carried out by induction on $k$ and $k'$, but we omit it.

For each vertex of the diagram $F_{{monster}}$, there is a corresponding vertex in the diagram $F'$. We did not describe this correspondence explicitly for every vertex, but it is easy to recover from the construction. We assign to the vertices of $F_{{monster}}$ the colors from $\mathcal B$ according to this correspondence.

To prove Proposition \ref{prop:OM}, it remains only to verify one property of the diagram $F_{monster}$, namely, that an even number of its vertices of  carry a symbol that projects to $2$.
These are precisely the vertices that are locally arranged as the vertex $A$ of $F'$. Let us count them.
\begin{itemize}
    \item one at $A$;
    \item one at $A'$;
    \item $k$ in the insertions of small copies of $L'$ in the right path;
    \item $k-1$ in the insertions in the right return path;
    \item $k'$ in the left path;
    \item $k'-1$ in the left return path.
\end{itemize}
The total number is $2k+2k'$. This number is even, and this completes the proofs of Proposition~\ref{prop:OM} and Theorem~\ref{thm:mute_odd_not_sofic}.

\section*{Acknowledgments} I.M. gratefully acknowledges support from the SNSF grant TMAG-2\_216487/1.

In addition, the authors would like to thank Egor Veprev, Pierre Guillon, and all participants of the S\'eminaire Rauzy for helpful discussions.

\end{document}